\documentclass[11pt,]{article}

\newcommand{\cor}{\vspace{.25in} \noindent {\bf Corollary.} \quad }
\newcommand{\df}{\vspace{.25in} \noindent {\bf Definition.} \quad}

\newcommand{\pf}{\vs \noindent {\it Proof:} \quad}

\newcommand{\rk}{\vspace{.25in} \noindent {\bf Remark.} \quad }
\newcommand{\thm}{\vspace{.25in} \noindent {\bf Theorem.} \quad }

\newcommand{\Z}{\mbox{\msbm{Z}}}

\newcommand{\be}{\begin{equation}}
\newcommand{\ee}{\end{equation}}
\newcommand{\ba}{\begin{array}}
\newcommand{\ea}{\end{array}}
\newcommand{\bea}{\begin{eqnarray}}
\newcommand{\eea}{\end{eqnarray}}
\newcommand{\bean}{\begin{eqnarray*}}
\newcommand{\eean}{\end{eqnarray*}}

\newcommand{\vs}{\vspace{.25in}}

\font\msbm=msbm10

\usepackage{amsmath}
\usepackage{amsthm}
\usepackage{amssymb}
\usepackage{epstopdf}
\usepackage{graphicx}
\usepackage{lscape}
\usepackage{longtable,rotating} 
\usepackage{pifont}
\usepackage{rotating}

\numberwithin{equation}{section} 

\DeclareRobustCommand{\gobblefour}[4]{}

\makeindex

\begin{document}

\title{
On Curves and Surfaces of Constant Width
}
\author{H. L. Resnikoff\footnote{Resnikoff Innovations LLC; howard@resnikoff.com}
}
\date{{\small 2013 May 31}}
\maketitle

\begin{abstract}
This paper focuses on curves and surfaces of constant width, with some additional results about general ovals. We emphasize the use of Fourier series to derive properties, some of which are known.
 
Amongst other results, we show that the  perimeter of an oval is $\pi$ times its average width, and provide a bound for the radius of curvature of an oval that depends on the structure of the harmonics in its Fourier series.

We prove that the density of a certain packing of Reuleaux curved triangles in the plane is 
$$\frac{2 (\pi - \sqrt{3})}{\sqrt{15} + \sqrt{7} - 2 \sqrt{3}} \simeq 0.92288$$
which exceeds the maximum density  for circles ($\simeq 0.9060$), and  conjecture this is the maximum for any curve of constant width.

For surfaces of constant width we show that  $ \rho_0(P)+ \rho_1(Q)=w$ where  the $\rho_i$ are the  principal curvatures, $P$ and $Q$ are opposite points, and  $w>0$ is the width. Moreover, an ovoid is a surface of constant width $w>0$ if and only if $\rho_{{\rm mean}}(P)+\rho_{{\rm mean}}(Q) =  w$ where $\rho_{{\rm mean}}(P)$ is the average radius of curvature at point $P$.

Finally, we provide a Fourier series-based construction that produces arbitrarily many new surfaces of constant width.
\end{abstract}


\section{Historical introduction}

Consider a smooth compact convex plane curve $\Gamma$, that is, a bounded convex curve in the euclidean plane whose cartesian coordinates $(x,y)$ are given by twice differentiable  functions of some convenient parameter. Without loss of generality assume that the origin $(0,0)$ lies inside the curve.

Such a curve is said to have {\it constant width}\index{constant width} if the distance between every pair of parallel  tangent lines is constant. This constant is called the {\it width} of the curve, denoted $w$. The simplest example is the circle of diameter $w$.

The simplest curves of constant width that are not circles are the family of curves first explicitly defined and studied by  \index{Reuleaux, Franz (1829-1905)}Franz Reuleaux, a German professor of mechanical engineering whose father and grandfather were builders of machines.   The constant width curves named for him were described in section 22, page 116 of  his classic work on the theory of machines  \cite{reuleaux 1876}, and their properties developed in section  25, page 129,  titled `Figures of Constant Breadth', and the following sections through section 29. A physical example is   shown in Model B02 of the \index{Reuleaux Model Collection}Reuleaux Model Collection at Cornell University.
 
The simplest Reuleaux curve is  the boundary of the intersection of three circles whose centers are the vertices of an equilateral triangle and whose radii are  equal to the side of the triangle. This `curved triangle', as Reuleaux called it,  must have been known in antiquity because it is so easy to construct -- it was probably known to the Greeks --  but it is less likely that its significance as a potential alternative  to the wheel was understood.

The thirteenth  century cathedral of Notre Dame in \index{Bruges}Bruges has Reuleaux triangle windows. It is said that triangles with curved edges appear in the notes of Leibniz and Leonardo da Vinci (I have not seen the documents), but that is a far cry from an explicit  statement and proof of their mathematical properties. 

The simplest curves of constant width were not the ones first studied by mathematicians.  This may seem curious, but ease of construction does not necessarily open the door to deeper mathematical properties, nor does it speak to the connection with applications. \index{Euler, Leonard (1707-1783)}Leonard Euler \cite{euler} was probably the first to study mathematical properties of non-circular curves of constant width, which he called \index{orbiform} {\it orbiforms}.  This work, his 513$^{{\rm th}}$ publication, was presented 235 years ago, in 1778. It was motivated by a problem from geometrical optics -- the construction of \index{catoptrix}{\it catoptrices}, closed convex curves with the property that a ray proceeding from an interior point returns to the point after two reflections from the curve.\footnote{The ellipse is a catoptrix.} Euler's investigation began with `triangular curves', about which we shall have more to say below. He proved that the evolvent of a triangular curve is a curve of constant width, and then used the curve of constant width to construct a catoptrix. The history is summarized in \cite{white}. Today, the catoptrix and triangular curves have been largely forgotten but the existence, if not the detailed properties of curves of constant width and their generalizations to higher dimension, are familiar to many mathematicians.

There are two theorems about curves of constant width that are of initial interest. One is universal:  the perimeter of a curve of constant width $w$ is equal to $\pi w$. Isn't that pretty! The other theorem is specific. It asserts that the Reuleaux curved triangle has the smallest area among the curves of a given constant width, whereas the well-known isoperimetric inequality shows that the circle has the largest area. The first theorem was proved by Barbier \cite{barbier} in 1860; the other, by Lebesgue \cite{lebesgue1} and, independently, by Blaschke \cite{blaschke} in 1915.

There have been interesting applications of curves of constant width.  A popular puzzle question asks why manhole covers are round. The answer is that a cover in the shape of a disk cannot fall through the hole. Nor can a cover whose perimeter has the shape of a curve of constant width. The  water department of the city of San Francisco distinguishes access ports in the roadway  for potable water from those for  reclaimed water by using circular  covers for the former and  covers in the shape of a Reuleaux triangle for the latter \cite{sfwd}. 

Curves of constant width can be used as profiles of cams that convert rotary motion into linear motion.  As a Reuleaux polygon $\Gamma$ rolls along a line $L$,   the point of contact with the line moves along $\Gamma$ until it reaches a vertex. The vertex stops moving while the polygon rotates around it (raising the centroid as it moves) until the next circular arc of the curved polygon is tangent to $L$. For instance, if the parameter is conceived as measuring time,  each vertex of the Reuleaux  triangle spends one-sixth of a revolution stationary, so that only one half of a revolution corresponds to forward motion of the point of contact. Thus a Reuleaux  triangle could be used, for example, to advance sprocketed motion picture film through a projector frame by frame, stopping each frame for a while (typically 1/30 second) so that the projector can flash it on the viewing screen. With the advent of digital projectors, this ingenious mathematics-based mechanical solution has become obsolete. 

A cunning application of curves of constant width is to coinage. The British 20p and 50p coins are 7-sided Reuleaux polygons, as is the Jamaican dollar and the 50 Fil coin from the United Arab Emirates. The Canadian dollar coin -- the `loonie' -- is an 11-sided Reuleaux polygon. These coins are easily distinguished by touch from circular ones of similar diameter, but like the latter they can pass through vending machines that measure diameter and weight to test validity.  But let us note that the Reuleaux curved triangle would have an additional advantage as a coin because, having the least area of all curves of the same constant width, it requires the least material and is least likely to wear a hole in your pocket.

There also several applications that have been mistakenly cited as examples of curves of constant width. Two examples:  the rotor of the Wankel rotary motor is similar to, but not, a curve of constant width \cite{wankel}, although a curved triangle of constant width could be used; and the Watts drill that, employing a bit in the shape of a Reuleaux triangle, drills holes that are very nearly but not perfectly  square (the corners are slightly rounded; cp. \cite{cox+wagon} for a bit that can drill a perfectly square hole).
\begin{center}
***
\end{center}
With regard to curves of constant width, little is new in this paper, but I would like to believe that its organization and the few original results may make it of some interest and use to the reader. Our point of view is similar to Fisher \cite{fisher}, and some results are necessarily the same. When similar methods are applied to the general oval they yield results that appear to be novel. 

The subject has attracted the interest of many distinguished mathematicians, and has a long and elaborate history, with papers still appearing in recent years. It is therefore difficult to say what might not be, somewhere, in the literature. As far as I know, new are the theorems about ovals in section \ref{sec: oval convex & curvature}; the material in section \ref{plane packing} on packing curves of constant width in the plane; the theorems about ovoids and surfaces of constant width in section \ref{surface curvature thm} and the method and new examples of surfaces of constant width in section \ref{surface examples}. Various formulae scattered throughout the text, of lesser interest,  may also be new.

The author is indebted to R. O. Wells for valuable comments on an earlier draft of this paper, and to Jayant Shah whose careful reading of several drafts of the paper led to significant improvements.

\section{Representation of an oval}
\subsection{Support function}
\label{support p}

The motivation for this paper was the naturalness with which \index{Fourier series} Fourier series can be used to study properties of curves and solids of constant width, and to produce examples of them. It became apparent that this method is more general, and can be applied to some extent to arbitrary ovals. With this change of context, many properties of curves of constant width appear as specializations of results about ovals. Our point of view is similar to Fisher \cite{fisher}, who observes that the use of Fourier series in this context goes back at least to Hurwitz \cite{hurwitz}.

For our purposes an \index{oval} {\it oval} $\Gamma$ is a simple closed plane curve -- $\Gamma$ does not intersect itself -- that is convex. Being a `curve', an oval is continuous. Typically ovals are egg-shaped curves -- the word derives from the Latin {\it ovus} for `egg'\footnote{In German, an oval is an `eilinie' -- an egg (shaped) curve.}  -- but, although the concept is an old one in mathematics it does not have a precise definition. So-called `Cartesian ovals', generalizations of ellipses,  were introduced by Ren\'{e} Descartes in 1637 in his studies of optics and later studied by Newton.  The ovals of Cassini, introduced in 1680, are quartic curves such that the product of the distances from two fixed points is constant but they are not necessarily convex.    Ovals are sometimes considered to have an axis of symmetry. Often the tangent is required to be a continuous function of position on the curve, so the curve does not have corners, but we shall not always insist on that requirement. According to our definition, a regular polygon is an oval, although this may lead to analytical complications. 

Suppose that $\Gamma$ is an oval, $P=(x,y)$  a point on $\Gamma$, and  $T_P$  the directed tangent line to $\Gamma$ at $P$.

Select a cartesian coordinate system whose origin $O$ lies in the interior of the oval. Let the perpendicular from $O$ meet the tangent line at $Q$ and let the line segment $\overline{O Q}$ make an angle $t$ with the positive $x$-axis.  Denote the distance from $O$ to $Q$  by $p(t)$. $p(t)$ is called the support function of the oval. The geometry is shown in  figure  \ref{geometry} on page \pageref {geometry}.

The equation of the tangent line can be considered as a function of the parameter $t$ and therefore the coordinates $(x,y)$ are also functions of $t$.  The two conditions: (1)\,$P$ lies on the tangent line, and (2)\,$\Gamma$ is the envelope of the 1-parameter family of tangent lines, enable one to solve for $x(t)$ and $y(t)$ in terms of the support function. 

Indicate differentiation with respect to $t$ by a prime. Then the parametric equations for $\Gamma$ are
\bea
x(t)&=& p(t) \cos t - p'(t) \sin t  \nonumber \\
y(t)&=& p(t) \sin t + p'(t) \cos t 
\label{xy equations}
\eea
These relations hold for all points of an oval where the tangent line is well defined. 

On one hand we have the general oval for which the support function varies from point to point, and on the other the circle, for which  it is constant and (if the origin is the center of the circle) equal to the circle's radius. The curve of constant width lies in between these extremes. This arrangement is a special case of a general classification of ovals according to decreasing symmetry properties of the support function, which will be developed in the next section.

\subsection{Fourier expansion}

The angular parameter $t$ is defined on an interval of length $2\pi$. If the function $p(t)$ is twice differentiable and periodic, it has a Fourier series expansion that converges to $p(t)$ at each point of the interval.  Weaker assumptions are also sufficient for the existence, convergence and differentiability  of the Fourier series. Write
\be
p(t) = a_0 + \sum_{k=1}^{\infty} a_k \cos  k t + b_k \sin k t 
\label{p fourier series}
\ee

Suppose $\Gamma$ is convex and consider all parallel lines that  have the same direction. Some will intersect $\Gamma$ twice and others not at all, but two of them, say $T_P$ and $T_{P^*}$, will be tangent to the oval. 

The {\it width} of $\Gamma$ at $P$ is the distance between $T_P$ and $T_{P^*}$ . In general, the width $w(t)$ depends on the point $P$, and hence on the parameter $t$; we shall call $w(t)$ the {\it instantaneous width} at $t$. The instantaneous width is expressed in terms of the support function $t \rightarrow p(t)$ by
\be
\label{p eq}
w(t) = p(t) + p(t + \pi)  
\ee
 If $w(t)$ is  independent of $P \in \Gamma$, then $\Gamma$  is said to be  a curve of constant width $w$.  

Denote the average width of the oval by $\overline{w}$.
\bea
\label{oval wbar}
\overline{w} &=& \frac{1}{\pi} \int_0^{\pi} w(t)\, dt \nonumber \\
&=&  \frac{1}{\pi} \int_0^{\pi} p(t) + p(t+\pi) \, dt \nonumber \\
&=&  \frac{2}{2 \pi} \int_0^{2 \pi} p(t) \, dt 
\eea
so the constant term in the Fourier series (which is also the average value $\overline{p}$ of the support function)  is $\overline{w}/2$:
\be
\label{oval FS}
p(t) = \frac{\overline{w}}{2} + \sum_{k=1}^{\infty} \big( a_k \cos k t + b_k \sin k t \big)
\ee
It follows that the Fourier series for the instantaneous width is
\be
\label{instantaneous width}
w(t) = \overline{w} + 2 \sum_ {k=1} ^ {\infty}    \big( a_{2k} \cos 2 k t + b_{2k} \sin 2k t \big)    
\ee
This series  expresses the instantaneous width in terms of the variation from the average width. Notice that only even harmonics appear in the expression for the instantaneous width. If the coefficients of the even harmonics are zero: $a_{2k}= b_{2k}=0$ for $k>0$, then the instantaneous width is equal to the average width  $\overline{w}$, that is, it is constant.  If the oval is convex this condition  produces a curve of constant width, and for the curve of constant width $w$ the support function takes the form
\be
\label{eq p gen sol}
p(t)=\frac{w}{2} + \sum_{k=0}^{\infty} a_{2k+1} \cos (2k+1) t + b_{2k+1} \sin (2k+1) t
\ee

\begin{center}
***
\end{center}

Every positive integer can be written in the form $2^n q$ where $q>0$ is odd. Hence the series eq(\ref {p fourier series}) can be re-written as
\be
p(t) = a_0 + \sum_{m=0}^{\infty} \sum_{q }a_{2^m q} \cos  \big(2^m q  t \big) + b_ {2^m q} \sin \big( 2^m q  t  \big)
\ee
This suggests an analysis based on the maximum power of 2 that occurs in the expansion of the support function.  

The $2^n$-th roots of unity are $ \big\{ \exp( 2 \pi i k/2^n) : 0 \leq k < 2^n \big\}$.  Consider ovals for which the quantity 
\be
\label{oval m-ave}
\overline{p}_m(t) := \frac{1}{2^m} \sum_{k=0}^{2^m-1} p(t + 2 \pi k/2^m)
\ee
does not depend on $t$.  This is the average value of the instantaneous support function when $t$ is augmented by the arguments of the $2^m$-th roots of unity. If an oval satisfies this condition for $m$, then it satisfies it for all $n>m$ because the condition for $n$ can be rewritten as $2^{n-m}$ copies of the condition for $m$. 

\thm
\label{thm oval radius}  If $\overline{p}_m(t) $ is constant, then $\overline{p}_m(t) =\frac{\overline{w}}{2}$, that is,
\be
\label{oval m mean p}
 \sum_{k=0}^{2^m-1} p(t + 2 \pi k/2^m) = 2^{m-1} \overline{w}
\ee
\pf Integrate the equation
$$ c = \frac{1}{2^m} \sum_{k=0}^{2^m-1} p(t + 2 \pi k/2^m) $$
to find
\bean c&=&  \frac{1}{2\pi} \int_0^{2 \pi}  \frac{1}{2^m} \sum_{k=0}^{2^m-1} p(t + 2 \pi k/2^m) \, dt\\
&=&  \frac{1}{2^m} \sum_{k=0}^{2^m-1}   \frac{1}{2\pi} \int_0^{2 \pi} p(t)\, dt \\
&=& \frac{\overline{w}}{2}
\eean $\Box$

\df The least $m$ for which $\overline{p}_m(t)$  is constant will be called the degree of the oval, and we also say $m$ is the degree of the Fourier series of the support function. An oval of degree $m$ is said to have constant width of degree $m$. $\Box$

If $m=0$, this condition just means that $\overline{p}_0(t)=p(t)=\overline{w}/2$ is constant so the oval is a circle. For $m=1$ it reduces to $\overline{p}_1(t)=p(t) + p(t+\pi)=\overline{w}$, a constant; these are the traditional curves of constant width. In the limit as $m \rightarrow \infty$ with $ \Delta t = 2 \pi /2^{m}$, the expression eq(\ref {oval m-ave}) approaches the Riemann sum 
$$ \lim_{m \rightarrow \infty} \frac{1}{2 \pi}\sum_{k=0}^{2^m-1}  p(t +k \Delta t ) \,  \Delta t  $$ 
for the integral
$$\frac{1}{2\pi}\int_0^{2\pi} p(t)\, dt $$
whose value we have seen is $\overline{w}/2$. 

Denote the set of ovals of degree $\leq m$ by $\mathfrak{O}_m$. Put $\mathfrak{O}_{-1} = \emptyset$.
Evidently
\be
\mathfrak{O}_m \subset \mathfrak{O}_{m+1}
\ee
The difference sets 
$$ \mathfrak{M}_{m} := \mathfrak{O}_{m}- \mathfrak{O}_{m-1}, \quad m \geq 0$$
are `pure' in the sense that the harmonics of the Fourier series (except for the constant term) are odd multiples of the same power of 2.

In terms of the previously considered Fourier series expansions,
$$ \mathfrak{M}_{0} = \big\{ p(t) = \overline{w}/2 \big\} , \quad \mbox{circle;}$$
$$ \mathfrak{M}_{1} = \Big\{ p(t) = \overline{w}/2 + \sum_{q \, \mbox{\tiny{odd}} }  a_q \cos q t + b_q \sin q t \Big\} , \quad \mbox{curves of constant width;}$$
$$ \mathfrak{M}_{2} = \Big\{ p(t) = \overline{w}/2 + \sum_{q \, \mbox{\tiny{odd}} }  a_q \cos 2q t + b_q \sin 2q t \Big\} , \quad \mbox{pure of degree 2;}$$
etc.

\thm The Fourier series for the support function of an oval has   degree $m$ fif\footnote{`fif' means `if and only if'.} it has the form
\be
\label{deg m}
p(t) = \frac{\overline{w}}{2} + \sum_{l=0}^{m-1} \sum_ {q \, \mbox{\tiny{odd}} } a_{l,q} \cos \big(2^l q t) + b_{l,q} \sin \big(2^l q t)
\ee
The pure Fourier series of degree $m$ has the form
\be
\label{pure deg m FS}
p(t) = \frac{\overline{w}}{2} +   \sum_ {q \, \mbox{\tiny{odd}} } a_{l,q} \cos \big(2^{m-1} q t) + b_{l,q} \sin \big(2^{m-1} q t)
\ee

\pf 
Suppose $p(t)$ has the given form. The constant term of 
$$\frac{1}{2^m} \sum_{k=0}^{2^m-1} p(t+ 2 \pi  k/2^m) $$
 is $\overline{w}/2$ and the sum over roots of unity for a typical summand is
$$ \sum_{k=0}^{2^m-1} a_{l,q} \cos 2^l q \big(t +  2 \pi k/2^l \big) +  b_ {l,q} \sin 2^l q \big(t +  2 \pi k/2^l \big) $$
Consider the cosine part of this expression:
$$ \sum_{k=0}^{2^m-1} \cos \Bigg( 2^l q t + \Big( \frac{2 \pi q k l}{2^{m-l}} \Big) \Bigg)
$$ 
The sum can be evaluated. The factors independent of $t$ are
$$ \frac{\sin 2^l q^2 \pi }{\sin 2^{l-m}  q^2 \pi } $$
The denominator is not zero because $l<m$. The numerator is zero because $2^l q^2 $  is odd. A similar argument applies to the sine part of the expression. Hence the series has degree $m$.

The proof of the converse is obtained by examining the general Fourier expansion, eq(\ref{oval FS}). Suppose that $\sum_{l=0}^{2^m-1} p(t + 2 \pi l/2^m) $ is constant, say $c$. Then

If $p(t)$ has degree $m$, then
\bean
c &=& \sum_{l=0}^{2^m-1} p(t + 2 \pi l/2^m)\\
&=&2^{m-1}  \overline{w}  + \sum_{l=0}^{2^m-1} \sum_{k=1}^{\infty} a_k  \cos k(t + 2 \pi l/2^m) +b_k  \sin k(t + 2 \pi l/2^m)  \\
&=&2^{m-1}  \overline{w}  + \sum_{k=1}^{\infty}  \Bigg( a_k  \sum_{l=0}^{2^m-1}  \cos k(t + 2 \pi l/2^m) +b_k  \sum_{l=0}^{2^m-1}  \sin k(t + 2 \pi l/2^m) \Bigg) \\
\eean
The inner sums are zero fif $k$ is an odd multiple of $2^m$ so the series reduces to a constant only when all other coefficients of the trigonometrical functions vanish. Indeed, as above, the argument reduces to consideration of $\sum_{l=0}^{2^m-1}  \cos  2 \pi k l/2^m$   and $\sum_{l=0}^{2^m-1}  \sin  2 \pi k l/2^m$, so consider
$$ \sum_{l=0}^{2^m-1} \exp  2 \pi i k l/2^m =  \sum_{l=0}^{2^m-1} \Big( u^{k/2^m}\Big)^l$$
where $k>0$. If $k/2^m$ is divisible by 2, the expression in parentheses is not a primitive root of unity. \\
$\Box$

The proof shows that the constant is $2^{m-1}  \overline{w}  $.
The simplest example of an oval of degree 2 is the lozenge-shaped curve shown in figure \ref{fig: deg 2 lozenge} on page \pageref {fig: deg 2 lozenge}, whose support function is
\be
\label{deg 2 lozenge}
p(t)=1 + \frac{1}{3} \cos 2t
\ee
The width varies from 4/3 to 8/3, with an average of 2, but the sum of the instantaneous width for any $t$ and the orthogonal instantaneous  width (for $t+\pi/2$) is constant, equal to 4.

\subsubsection{Convexity and curvature}
\label{sec: oval convex & curvature}

$\Gamma$ will be convex fif its curvature does not change sign. In this case, we can arrange that the curvature is $\geq 0$ without loss of generality.

The curvature $\kappa$  is given by
$$
\kappa(t) = \frac{x'(t) y''(t) - y'(t) x''(t)}{\big(x'(t)^2 + y'(t)^2 \big)^{3/2}}
$$
and the radius of curvature $\rho(t)$ is its reciprocal: 
\be
\label{rho1}
\rho(t) = 1/\kappa(t)
\ee

In the present situation -- because the fundamental expression is the support function rather than cartesian or polar coordinates --  these formulae simplify to attractive expressions. One finds
$$
\label{rho2}
\rho(t) = \sqrt{ \big( p(t) + p''(t) \big)^2 }
$$
Since we are only interested in convex  curves, for which the expression is not negative, the square root can be extracted:  then $\rho(t) = p(t) + p''(t)$. 
If the curve is closed the condition 
\be
\label{p convexity}
\rho(t) = p(t) + p''(t) \geq 0 
\ee
 is equivalent to convexity of $\Gamma$. 
 
Notice that if  $t$ is thought of as representing time, the support function $p(t)$ can be thought of an harmonic oscillator driven by the radius of curvature $\rho(t)$.\footnote{This leads to a method for constructing ovals with large maximum radius of curvature which will not be elaborated here.}

Expressing $\rho(t)$ as a real Fourier series implies that $\Gamma$ is convex fif
\be
\rho(t) = \frac{\overline{w}}{2} + \sum_ {k=2} ^ {\infty} ( 1-k^2 )  ( a_k \cos k t + b_k \sin k t  ) \geq 0
\label{oval rho}
\ee
For instance, if the average width of the curve is positive and 
$$\sum_{k>1} ( k^2-1)\big( |a_k|^2 + |b_k|^2 \big) \leq  \frac{\overline{w}}{2}$$ then the curve is convex, that is, an oval.

There is another connection between the radius of curvature and the width of an oval:

\thm Suppose that an oval of degree $m$ can be expressed in terms of a  support function $p$ that has a convergent and twice differentiable Fourier series. Denote the average width of $\Gamma$ by $\overline{w}$ and  the average  radius of curvature by $\overline{\rho}$. Then
\be
\label{oval m rho root sum}
\sum_{k=0}^{2^m-1} \rho\big( t + 2 \pi k/2^m \big) 
=  2^{m-1}  \overline{w}  
\ee
\be
\label{oval m rho mean}
\overline{\rho} = \frac{\overline{w}}{2}  
\ee
\be
\label{oval m rho bound}
0  \leq  \rho(t) \leq   2^{m-1}  \overline{w}  
\ee

\pf 
Equation (\ref {oval m mean p}) states that if the degree of $p(t)$ is $m$  then
$$ \sum_{k=0}^{2^m-1} p\big( t + 2 \pi k/2^m \big) = 2^{m-1} \overline{w}$$
whence
\bean
 \sum_{k=0}^{2^m-1} \rho\big( t + 2 \pi k/2^m \big) 
&=& 
  \sum_{k=0}^{2^m-1} p\big( t + 2 \pi k/2^m \big)+  \frac{d^2}{d\,t^2} \sum_{k=0}^{2^m-1} p\big( t + 2 \pi k/2^m \big)\\\
  &=&  2^{m-1} \overline{w}\\
\eean
which is the first result. 

The average  radius of  curvature is
$  \overline{\rho}=\frac{1}{2 \pi}  \int_0^{2\pi}  \rho(t) dt  $.  Integrate to find:
\bean
 2^{m-1} \overline{w} & =& \frac{1}{2\pi} \int_0^{2\pi}  \sum_{k=0}^{2^m-1} \rho\big( t + 2 \pi k/2^m \big) dt   \\
 &=& \sum_{k=0}^{2^m-1}  \frac{1}{2\pi} \int_0^{2\pi} \rho\big( t + 2 \pi k/2^m \big) dt \\
 &=&   \sum_{k=0}^{2^m-1}  \frac{1}{2\pi} \int_0^{2\pi} \rho( t) dt \\
 &=& 2^m \overline{\rho}
\eean
which completes the proof of the second assertion.

For the last part of the theorem, recall that the radius of curvature of an oval is non-negative. If all but the first of the terms in eq(\ref {oval m rho root sum}) is zero, the result follows.  \\
$\Box$

\cor The average of $\rho$ with respect to the roots of unity  is equal to the average of $\rho$ over the circle, i.e.
\be
\overline{\rho}_m = \overline{\rho}
\ee

\rk If the oval has symmetries, the maximum radius of curvature will be repeated some number of times, which will tighten the bound on it. For instance, the lozenge shown in figure \ref{fig: deg 2 lozenge} has degree $m=2$ (the support function is $p(t)=1+\frac{1}{3} \cos 2t $). According to the theorem, $\rho(t) \leq 4$, but the bilateral symmetry  implies that the maximum radius of curvature occurs twice so $\rho(t) \leq 2$. Indeed, the maximum, 2, and the minimum, 0, are both achieved.

\rk  The theorem applies in particular to curves of constant width $w$ for which $m=1$ and $w=\overline{w}$.  In this case the third part states that $0 \leq \rho(t) \leq w$. How can we reconcile this with the fact, evident from its construction, that the Reuleaux curved triangle of width 2 consists of three congruent circular arcs of radius 2 that meet at three points --  the `vertices' of the curved triangle? 

The answer is that the radius of curvature is zero at the vertices, and as a point traverses the curve according to the parametric equations, it dwells at the vertices for a `time' equal to the time spent traversing the circular arcs. $\Box$

The final assertion of this theorem, eq(\ref {oval m rho bound}), is the key for building ovals that have a large radius of curvature: one must employ high degree. For ovals such as polygons,\footnote{in this case, the radius of curvature has an infinite discontinuity.}  that have flat spots where the radius of curvature is infinite, the degree must be infinite.
 \begin{center}
\rule{60pt}{0.5pt}
  \end{center}
 
\begin{small}
The reader may be interested in the details of the construction of ovals with large radius of curvature. For polygons, the radius of curvature at each point of the straight edges is infinite but in terms of the support function, each edge corresponds to a single parameter value and the radius of curvature is realized by a distribution supported there (a `Dirac delta' distribution). The support function spends all of its time `turning the corners' of the polygon's vertices. For ovals that have a finite but arbitrarily large radius of curvature at some point traversal of the curve by the parameter `speeds up' as the radius of curvature increases and `slows down' as  the curvature decreases.  

In the construction of the oval a difficulty can arise because the Fourier series terms oscillate around their limit and finite approximations could, therefore, cause either the support function or the curvature to be negative on some intervals. Thus the oval would be approximated by a sequence of curves that approach ovality but are not themselves ovals. This would require consideration of infinite Fourier series with the attendant analytical complexity and issues of convergence.

However, it is possible to construct ovals with arbitrarily large radius of curvature whose support function is a finite Fourier series. In the study of convergence of Fourier series  the Fej\'{e}r kernels $F_n(t)$ play a role. Let $n>0$ be a an integer and define
\be
\label{Fejer1}
F_n(t) = \frac{1}{n} \Big( \frac{1 - \cos (n t) }{1 - \cos t} \Big)
\ee
The Fej\'{e}r kernel is non-negative. 

For our purposes the expression of the Fej\'{e}r kernel as a finite Fourier series whose maximum lies at $t=\pi$  is preferable:
\be
\label{Fejer}
F_n(t) = 1 + 2\sum_{k=1}^{n-1} (-1)^k \frac{n-k}{n} \cos( k t)
\ee
The maximum of $F_n(t)$ is $n$ and its average value is 
$$\frac{1}{2 \pi} \int_{0}^{2\pi} F_n(t) = 1 $$
Its graph consist of one tall hump centered at $t=\pi$ surrounded by many small oscillations (there are $n-1$ minima, all non-negative). Thus we may try to construct an oval using this function to represent the curvature.

There is a problem. The Fourier series for the  Fej\'{e}r kernels contains a term proportional to $\cos t$ which cannot arise from a finite Fourier series $p(t)$. This term can be avoided by introducing symmetry into the oval. Let $\sigma > 1$ denote an integer. The graph of $F_n(\sigma t)$ will have $\sigma$ congruent peaks. Now the support function can be found either by solving the differential equation
$$ p(t)+p''(t)= F_n(\sigma t), \quad p(0)=p(2\pi) $$
or by observing that the Fourier series for $p(t)$ can be directly written down. It is
\be
p(t)= 1 + 2 \sum_{k=2}^{n-1} (-1)^{k}  \Big(\frac{n-k}{n} \Big) \Big(\frac{1}{1-(\sigma k)^2} \Big) \cos( \sigma k t)
\label{Fejer p n s}
\ee
That $p(t)$ is non-negative follows from examining the majorant 
$$2 \sum_{k=2}^{n-1}  \Big(\frac{n-k}{n} \Big) \Big(\frac{1}{1-(\sigma k)^2} \Big) $$
of the sum, which can never be less than -1 if $\sigma \geq 2$.

The oval defined by the support function of eq(\ref {Fejer p n s} ) has maximum radius of curvature $n$ and  its shape is a curved  $\sigma$-sided  regular polygon that has the regular $\sigma$-sided polygon as limit for $n \rightarrow \infty$. For instance, figure \ref{fig: Fejer n=32,sigma=5} on page \pageref{fig: Fejer n=32,sigma=5} shows a pentagonal oval ($\sigma = 5$) whose maximum radius of curvature is $n=32$.

\end{small}
 \begin{center}
\rule{60pt}{0.5pt}
  \end{center}
 
\subsubsection{Arc length}

Denote the arc length parameter by $s(t)$. The perimeter of the oval $\Gamma$ is
\be
s = \int_{\Gamma} ds = \int_{0}^{2\pi} \sqrt{x'^2 + y'^2} \, d t
\ee
From eq(\ref{xy equations}) the derivatives are

\bea
x'(t) &=& - (p(t)+  p''(t))  \sin t  \nonumber  \\
y'(t) &=& - (p(t)+  p''(t))  \cos t  
\label{parametric eqs}
\eea
whence
\be
x'^2 + y'^2 =\big( p(t)+  p''(t) \big)^2
\label{p arc length element}
\ee
If $p(t)+  p''(t) \geq 0$, that is, if $\Gamma$ is convex, then the integrand simplifies to $p(t)+  p''(t)$ and then, subject only to the necessary conditions for absolute convergence of the series should it have an infinite number of terms, it follows that 
\bea
s &=& \int_{0}^{2\pi} p(t)+  p''(t)   \, d t \nonumber  \\
 &=&  p'(2\pi)-p'(0)  + \int_{0}^{2\pi} p(t) \, d t   \nonumber \\
  &=& \int_{0}^{2\pi} p(t) \, d t \nonumber \\
   &=&  \pi \overline{w} + \sum_{k=1}^{\infty} \Big( a_k  \int_{0}^{2\pi}  \cos k t \, d t + b_k  \int_{0}^{2\pi}  \sin  k t \, d t  \Big) \nonumber \\
    &=&  \pi \overline{w}
\eea
This proves
\vspace{-10pt}
  
 \thm The perimeter of an oval is $\pi$  times  its average width and $2 \pi$ times its average radius of curvature.
 \label{oval Barbier thm} 
 \vspace{6pt}

Thus all ovals that have the same average width have the same perimeter. This is a generalization of Barbier's well-known  theorem for curves of constant width, which asserts that all curves of given constant width $w$  have perimeter $\pi w$. Barbier's proof \cite{barbier} of his pretty theorem, and some by later mathematicians, employ probabilistic arguments. The analogue of the theorem of Barbier  is false for surfaces of constant width in 3-space: the surface area of a solid of given constant width is not independent of the volume of the solid.
 \begin{center}
\rule{60pt}{0.5pt}
  \end{center}
  
 \begin{small}
 An example of the theorem for an oval whose width is not constant and whose support function has a Fourier series whose harmonics have arbitrarily large powers of 2 as factors may be of interest. Consider the ellipse given by
 $$ \frac{x^2}{a^2} + \frac{y^2}{b^2} =1, \quad a>0, \, b>0 $$
 The support function relative to the origin of coordinates is
 $$ p(t) = \Big( \big(a \cos t \big)^2 + \big(b \sin t \big)^2  \Big)^{1/2} $$
 Direct calculation of the average width yields\footnote{We use {\it Mathematica}\texttrademark to calculate integrals.}
 \bea
\pi  \overline{w} &=& \frac{1}{2} \int_0^{2\pi} p(t)+p(t+\pi) \,  dt  \nonumber \\
  &=& 2  \Bigg( a  \,E \Big( 1 - \frac{b^2}{a^2}\Big) + b  \, E\Big( 1 - \frac{a^2}{b^2}\Big) \Bigg)
  \label{width ellipticE}
 \eea
 where $E(m)=\int_0^{\pi/2} \sqrt{1 - m \sin^2 \theta} \, d\theta $
  denotes the complete elliptic integral of the second kind, while calculation of the perimeter $s$ using the parametrization eq(\ref{xy equations}) of $(x,y)$ provided by the support function  produces
\be
s = 4 a E \Big( 1 - \frac{b^2}{a^2} \Big) 
\label{arc ellipticE}
\ee
Although they appear to be different,  the expressions eq(\ref{width ellipticE}) and eq(\ref{arc ellipticE}) for $\pi \overline{w}$ and the perimeter $s$ are the same because of the known identity
$$  E\Big(1-z^2 \Big) = z E\Big(1-z^{-2} \Big) $$
 \end{small}

\subsubsection{Area}

The area bounded by an oval can be calculated directly from the parametric Fourier series representation of the curve.

 If the curve is presented by differentiable parametric equations $x=x(t), \, y=y(t)$, $0 \leq t <2 \pi$,  then the area it bounds is given by the well-known formula for the area of a region bounded by a convex closed curve:
\be
\label{area-xy}
A_{\Gamma} = \frac{1}{2} \int_{0}^{2\pi}  r^2  d \theta
\ee
where $(r,\theta)$ are the polar coordinates of ($x,y)$. 
In the present case $x=x(t)$ and $y=y(t)$ are {\it not} given as function of the polar angle $\theta$ so we must find the relationship between $\theta$ and $t$. A diagram shows $ \tan (t-\theta)  = -p'(t)/p(t) $ whence $ \theta = t + \arctan (p'/p) $
and 
$$d\theta = \Bigg( 1+\frac{p(t) p''(t) - p'(t)^2}{p(t)^2+p'(t)^2} \Bigg) dt$$
Noting that $r^2 = x(t)^2+y(t)^2=p(t)^2+p'(t)^2$, 
$$
A_{\Gamma} = \frac{1}{2} \int_0^{2\pi} p(t)^2 +  p'(t)^2 dt + \frac{1}{2}
 \int_0^{2\pi} p(t) p''(t) - p'(t)^2 dt
$$
In the second integral the first term can be integrated by parts. Noting that the endpoints make no contribution,  the result is
\be
\label{oval area}
A_{\Gamma} = \frac{1}{2} \int_0^{2\pi} p(t)^2 -  p'(t)^2 dt 
\ee

We could insert the Fourier series and calculate the integrals in the double infinite series. But this is a good moment to recall that  the Fourier basis functions $\{\frac{1}{\sqrt{2 \pi}} \} \cup \{ \frac{1}{\sqrt{ \pi} }  \cos k t  : k>0 \}\cup \{ \frac{1}{\sqrt{ \pi} }  \sin k t : k>0 \}$ are a complete  orthonormal system of functions in the space of square integrable functions on $[0,2\pi)$, and to recognize that 
$ \int_{-\pi}^{\pi} p(t)^2 dt$ is the squared norm of $p(t)$ and $ \int_{-\pi}^{\pi} p'(t)^2 dt $ is the squared norm of $p'(t)$. The squared norm of a function  is just the sum of the squared coefficients of its Fourier series with respect to the orthonormal basis functions. So, assuming that the average width of the oval is positive, we can immediately write 
\be
\label{oval area}
A_{\Gamma}=   \frac{\pi}{4} \Big( \overline{w}^2 - 2 \sum_{k=2}^{\infty} \big( k^2 - 1 \big)\big(a_{k}^2 + b_{k}^2 \big)\Big)
\ee
 This agrees, as it must,  with the isoperimetric inequality, which implies that the area of the circle is greater than the area of any other oval of the same average width. 

\rk This formula leads to a curious result. Suppose that the support function is a Fourier series with a finite number of summands. If all the non-zero Fourier coefficients are rational, then the area of the oval is a rational multiple of $\pi$. The perimeter is also a rational multiple of $\pi$, namely:  $\overline{w} \pi$.  And the average radius of curvature is the  rational number $\overline{w}/2$.

\section{Curves of constant width}

According to eq(\ref {eq p gen sol}),  an oval is a curve of constant width $w$ if the Fourier series of its support function has the form
$$
p[(t) = \frac{w}{2} + \sum_{k=0}^{\infty} a_k \cos (2k+1) t + b_k \sin (2k+1) t
$$

\subsection{Convexity again}

Although the radius of curvature of an oval of given average width is not necessarily bounded,  one sees by specializing  eq(\ref{oval m rho bound}) to  degree $m=1$, that for a curve of constant width $w$,
\be
0 \leq \rho(t) \leq w
\ee

\subsection{Mellish's theorem}

According to the definition, an oval is a curve of constant width $w$ fif the distance between parallel tangent lines is $w$. Mellish \cite{mellish} gave a similar criterion in terms of the radius of curvature $\rho(t)$:  an oval is a curve of constant width fif 
\be
\label{Mellish thm}
\rho(t) +\rho(t+\pi) = w
\ee

From the general Fourier series we calculate
\be
\label{p radius of curvature}
\rho(t) = p(t) +  p''(t) 
\ee
Consider $\rho(t)+\rho(t+\pi)$.  Now $ p''(t) +p''(t + \pi) = 0$ because the derivative kills the constant term in the series and the trigonometric functions change sign when  $\pi$ is added to the argument. The surviving terms are $p(t) +p(t+\pi) =w$.  $\Box$

An interesting consequence of Mellish's theorem is that the maximum and minimum values of the radius of curvature sum to the width. Denote the maximum, resp. minimum,  of $\rho(t)$ by $\rho_{{\rm max}} $, resp. $\rho_{{\rm min}} $ and suppose for the moment that $\rho(t)$ is twice differentiable. Then \label{rhoMin+rhoMax}
\be
\label{radius max min}
w = \rho_{{\rm max}} +\rho_{{\rm min}}
\ee
The extrema of $\rho(t)$ are the roots of $\rho'(t)=0$. Suppose $t_0$ is such a root. Then
$$ \rho'(t_0)=\sum_k (2k+1)   \big( a_{2k+1} \cos (2k+1)t_0 + b_{2k+1} \sin (2k+1) t_0 \big) = 0 $$
so $\rho'(t_0+\pi)=0$ as well: $t_0+\pi$ is also an extremum. Testing the second derivative to determine its character, one sees that $ \rho''(t_0)$ and $ \rho''(t_0 + \pi)$ sum to zero: one is a local maximum, the other a local minimum. Consider a global maximum. It will be paired with some  minimum, but that minimum must be a global minimum for if not, the global minimum would be paired with a local maximum that is greater than the global maximum (because the sum of the pairs equals the width),  which is a contradiction.  It all depends on the absence of even harmonics in the Fourier series for $p(t)$.

\subsection{Examples}

\subsubsection{Examples of low degree}

The constant support function describes a circle centered at the origin. The fundamental frequencies of the Fourier series $\cos t$ and $\sin t$ translate the center of the circle. 

We can think of each of the higher harmonic terms of the Fourier series as contributing  an alteration of the circle that -- under some conditions that insure convexity of the result -- produces an oval. It the components of even order are avoided, the alteration will preserve the constant width of the circle.

The average width of $\Gamma$ depends only on the constant term. This implies that we must select the constant   $\overline{w}>0$ if we are to have an oval with positive width.  Later, in section \ref {sec: zero width},  we shall investigate what it means for an oval to have zero width.

For the rest of this section assume  without loss of generality  that the average width is $\overline{w}=2$.

The first terms that can change the shape of the circle are the second harmonics. Consider $a_2 \cos 2 t$. The graph of this component reminds one of an astroid -- a 4-cusped hypocycloid.  The parametric equations of the oval are
 \bean
 x(t) &=&(1 + 2a_2 -a_2 \cos 2t)\cos t   \\
 y(t) &=&(1 - 2a_2 -a_2 \sin 2t)  \sin t 
 \eean
 The curvature is $1-3 a_2 \cos^2 (2t) $ which is $\geq 0$ fif $|a_2| \leq 1/3$. These curves are all convex lune-like curves of average width 2 whose elongation varies from the circle to an extreme ratio of greatest to least width of 2:1. Each oval has perimeter $2 \pi$.
 
The simplest examples of  components that preserves constant width 
are the third harmonics. Consider the component  $ p(t)=a_3 \cos 3t$ with zero constant term. Application of trigonometric identities yields
\bea
x^c_3(t) &=& \hspace{7 pt} a_3(2 \cos 2t -  \cos 4 t )   \nonumber    \\
y^c_3(t) &=& -a_3(2 \sin 2  t + \sin 4 t  )
\eea
This curve looks like an equilateral triangle with concave sides. It is not smooth and it does not have a `width' as we have defined it; it is curve of `zero width' in the sense of Euler. The associated curve of constant width ($w=2$) for which $p(t)=1 + a_3 \cos 3t$ is  
\bean
x(t) &=&   \Big(  \cos t + 2 a_3  \cos 2t -   a_3 \cos 4 t \Big)    \\
y(t)   &=&   \Big( 1 - 12 a_3 \cos 2  t  -4 \cos 4 t  \Big) \sin t 
\eean

If  the average width is positive (in this case we have agreed to normalize it equal to 2) and the parameter $a_3$ is 0, then $\Gamma$ is a circle; we shall not be surprised to find that for sufficiently small values of $a_3$ the curve has constant width 2. We shall see that  the curve is convex fif $a_3 \leq 1/8$, and it is smooth if $a_3 < 1/8$. The radius of curvature is
$$p(t) + p''(t)=  1-8 a_3 \cos 3t $$
$
\Gamma$ is convex -- the curve is an oval --   fif $|a_3| \leq 1/8$.  Figure \ref{Rab3} displays this curve of constant width for the largest value of $a_3$  color-coded by hue on the color circle to show how the radius of curvature varies along the curve, from 0 at the three vertices to 2 at the midpoints of the three arcs. The circle of the same diameter is also shown; its radius of curvature is 1, indicated by cyan.

\subsubsection{Equation for Rabinowitz' curved triangle}

Rabinowitz \cite{rabinowitz} found a polynomial equation for a `curved triangle' of constant width  which is the special case  of $a=1/9$ of the family of curves $\Gamma_a$ defined by  $p(t) = 1+a \cos(3t)$. In this case  the curve has `rounded corners'; cp. eq(\ref{poly eq rabinowitz}).   $\Gamma_a$ is convex, and hence a curve of constant width $2$, whenever $a \leq 1/8$. In the limiting case, the `corners' have  infinite curvature. But $\Gamma_{1/8}$ is not Reuleaux's curved triangle, which has `corners' where the curvature is discontinuous but its limits on opposite sides of the corner are finite, like a corner of a square.

The area of the Rabinowitz curve of parameter $a$ is
\be
A(a) = (1-4 a^2) \pi
\ee
The area decreases from $\pi$ (corresponding to  $a=0$ and the circle of radius 1) to $\frac{15}{16}\pi$ for the limiting curve of constant width corresponding to $a =1/8$.  By the Blaschke-Lebesgue Theorem \cite{blaschke,lebesgue1,lebesgue2}, the Reuleaux curved triangle has the smallest area for a curve of its constant width, namely, $2( \pi -\sqrt{3}) \sim 2.819$ whereas $\frac{15}{16}\pi \sim 2.945$. 
 \begin{center}
\rule{60pt}{0.5pt}
  \end{center}

\begin{small}
Following the procedure in \cite{rabinowitz}, we find the following unpleasant polynomial equation of degree 8 for the curve of constant width defined by $p(t) = 1 + a  \cos 3 t$ (with $ |a| \leq 1/8$)  as a function of the Fourier coefficient parameter $a$:

\bea
\label{poly eq rabinowitz}
0 &=&-432 a^5 \left(x^3-3 x y^2\right)-243 a^6 \left(4 x^2+4 y^2+9\right)+ \nonumber  \\
	&& 27 a^4 \left(10 x^4+x^2 \left(20  y^2-9\right)+10 y^4-9 y^2+81\right)+\nonumber   \\
	&& 72  a^3 x \left(4 x^4-x^2 \left(8  y^2+15\right)-12 y^4+  45 y^2\right)+ \nonumber  \\ 
	&&2  a x \left(8 x^6-4 x^4 \left(2  y^2+9\right)+x^2 \left(-40 y^4+   72  y^2+27\right)- \right.  \nonumber \\
	&& \qquad \left. 3 y^2 \left(8 y^4-36   y^2+27\right) \right)+  \nonumber  \\
   	&& a^2 \left(100   x^6-3 x^4 \left(92 y^2+171\right)+9 x^2  \left(76 y^4-114 y^2+135\right)+ \right.  \nonumber \\
	&& \qquad  \left.   9
   \left(4 y^6-57 y^4+135  y^2-81\right)\right)+ \nonumber   \\
   	&&729 a^8+\left(x^2+y^2-1\right)  \left(x^2+y^2\right)^3
   \eea
   
If $a=0$ this reduces to
   \be
   0 = \left(x^2+y^2-1\right)
   \left(x^2+y^2\right)^3
   \ee
which is the equation of a circle of diameter 2 supplemented by the 6-tuple point $(0,0)$. The possible production of spurious factors is a concommitant of the method.

If $a=1/8$, which is the limiting case for an oval,  the resulting equation is

\bea
0 &=& 16777216 x^8+33554432 x^7+x^6 \left(67108864
   y^2+9437184\right)+  \nonumber \\
  &&  x^5 \left(-33554432
   y^2-141557760\right)+  \nonumber  \\
   && x^4 \left(100663296
   y^4-122683392 y^2-133373952\right)+   \nonumber \\
   &&x^3
   \left(-167772160 y^4+283115520
   y^2+77635584\right)+  \nonumber \\
   &&x^2 \left(67108864
   y^6+128974848 y^4-266747904 y^2+317447424\right)+  \nonumber \ \\
   &&x \left(-100663296 y^6+  
    424673280 y^4-232906752
   y^2\right)+  \nonumber  \\
   &&16777216 y^8-7340032 y^6-133373952
   y^4+317447424 y^2    \nonumber \\
   && -182284263 
\eea

Is this a polynomial equation for Reuleaux's curved triangle? Does $a=w/16$ yield his curve of constant width?  We already know the answer is `no' because the area is too large.  This is an equation for a curve of constant width, but not for the Reuleaux curved triangle. 
\end{small}

\subsubsection{Reuleaux's curved polygons of constant width}

Reuleaux's curved polygons of constant width are convex curves that consist of an odd number $q>1$ of congruent circular arcs which connect the vertices of a regular polygon of $q$ sides. The center of each circular arc is the vertex most distant from the  arc's endpoint vertices. Thus the curve has $q$-fold rotational symmetry, and it is also symmetric about each line that passes through the origin and a vertex.

For any curve of constant width, set the width $w=2$ and the  coefficients of $\cos t$ and $\sin t$  to zero. Then the Fourier series representation eq(\ref{eq p gen sol}) of the support function $p(t)$ will contain only cosine functions (because of the reflection symmetry) and will contain harmonics of $q$ because of the $q$-fold rotational symmetry. Thus the general term will be proportional to $\cos q k t$. Since only odd harmonics can appear, both $q$ and $k$ must be odd. Hence
\be
p(t) = 1 + \sum_{n=0}^{\infty} c_{n} \cos q(2n+1)t, \quad q \, \mbox{odd}
\ee
This applies in particular to the Reuleaux curved polygons.

Reuleaux's curved triangle is the simplest to construct geometrically but complicated to study  analytically. It does not satisfy our assumptions because the limits of the tangent at three points that are the vertices of an equilateral triangle are unequal. Between these vertices the curve is merely an arc of a circle whose radius is the side of the triangle, hence smooth. 

When the parametric equations of a curve are known, the easiest way to calculate the width function $p(t)$ is from eq(\ref{xy equations}). We find
\be
p(t) = x(t) \cos t +  y(t) \sin t
\ee
For the Reuleaux curved polygons it is easy to express  $x$ and $y$ explicitly, and from that, to find $p(t)$. The Fourier series for $p(t)$ for Reuleaux's curved triangle of width 2 is
\be
\label{reuleaux triangle}
p(t)= 1 + \left( 1- \frac{2}{3}\sqrt{3} \right) \cos t + \frac{1}{\pi} \sum_{n=0}^{\infty} \frac{(-1)^n}{(2 n+1) (3 n+1) (3 n+2)} \cos (3(2n+1)t
\ee
A closed-form expression can be given in terms of the hypergeometric function $_2F_1(a,b,c,z)$.

It follows from eq(\ref{oval area}) that the area of a Reuleaux triangle of width 2 is
\be
\label{area Reuleaux triangle}
A= \pi - \frac{2}{\pi} \sum_{n=0}^{\infty} \frac{1}{(2n+1)^2 (3n+1)(3n+2) }  = 2(\pi -\sqrt{3}),
\ee
a result that is much easier to obtain directly from the geometry.

Figure \ref{Reu3_rho(t)_polarplot} shows how the radius of curvature with the parameter $t$. It alternately assumes the values 0 and 2 as $t$ varies along intervals of length $\pi/6$. This means that the point $(x(t),y(t))$ dwells at the vertcies of the curve half the `time' and moves with constant `speed' the other half. As always, the sum of the radii of curvature at opposite points $\big((x(t), y(t)\big)$ and $\big(x(t+\pi),y(t+\pi)\big)$ is 2.

\subsubsection{Nowhere differentiable curvature}

Before the last decades of the 19$^{{\rm th}}$ century it was generally believed that a continuous function must be differentiable except for a limited number of `corners' where left- and right-derivatives have different limits. In 1872 Weierstrass constructed a family of  continuous functions that were nowhere differentiable.  The class of ovals whose radius of curvature is continuous but nowhere differentiable provides an instructive collection of examples. 

We shall limit our discussion to Weierstrass's functions of the form
\be
\label{W-function}
W_{a,b}(t) = \sum_{k=1}^{\infty} a^k \cos b^k t
\ee
where $ 0<a<1$, $b$ is an integer greater than 1, and $a b \geq 1$. In these circumstances (which are more general than Weierstrass's  conditions; see  \cite{Hardy}),   $W_{a,b}(t)$ is continuous and nowhere differentiable.

Consider an oval with support function $p(t)$ and, therefore, radius of curvature
$\rho = p + p''$. Suppose that
\be
\label{rho ND}
\rho(t) = 1 + c \,W_{a,b}(t)
\ee
The constant $c$ will be determined later. The graph of $\rho$ is fractal. The support function can be reconstructed from this expression; it is
\be
\label{p ND}
p(t) = 1 + c \sum_{k=1}^{\infty} \Big( \frac{a^k}{1-b^{2k}} \Big) \cos b^k t
\ee
We know that such a series produces a curve of constant width $w=2$ fif all the frequencies are odd, that is, fif $b$ is an odd integer. Suppose this is so. Then the maximum radius of curvature will be achieved for $t=0$ and it will have the value 
$$\rho_{{\rm max}} = 1 + \frac{a c}{1-a} $$ 
The minimum of the Weierstarss function will be achieved for $t=\pi/b$ because then $\cos (b^{k-1} \pi )=-1$; here is where we use the fact that the frequency is odd. The minimum  will be 0 fif
$$\rho_{{\rm min}} = 1- \frac{a c}{1-a} =0$$
so   $c= \frac{1-a}{a} $ whence the largest radius of curvature  is 
$$\rho_{{\rm max}} = 2$$
 in agreement with our theorem about the range of the radius of curvature of a curve of constant width $w=2$. Thus no matter what parameter $a$ and  odd integer $b$ are chosen, the corresponding Weierstrass curve of constant width realizes the possible extremes 0 and 2 for the radius of curvature.

But we do not have to restrict ourselves to the infinite Fourier series whose graph is a fractal. If the series is truncated after $N$ terms, the radius of curvature will be continuous and differentiable and a similar analysis yields curves of constant width whose maximum radius of curvature approaches 2 as $N \rightarrow \infty$.

What happens if 2 divides $b$? In this case we obtain an oval of average width $\overline{w}=2$ that is not of constant width. The minimum and maximum radii of curvature can be calculated in the same way as before. For instance, for $a=1/2, b=2$, we find the minimum occurs for $t=\pi/3$ and $c=2$, so
$$  \rho_{{\rm max}} = 3 , \quad  \rho_{{\rm min}} = 0 $$

\subsection{Curves of zero width}
\label{sec: zero width}

Let us reconsider eq(\ref{instantaneous width}), which expresses the {instantaneous width $w(t)$  of the oval $\Gamma$ as a sum of deviations, all of which are even harmonics,  from the average width $\overline{w}$:
\be
\label{instantaneous width2}
w(t)   =  \overline{w}  + 2 \sum_ {k=1} ^ {\infty}    \big( a_{2k} \cos 2 k t + b_{2k} \sin 2k t \big)    
\ee

If  $\overline{w}>0$ then we can renormalize it --  choose a unit of measure so that   $\overline{w}=2$ -- and this is the situation we have explored.   For $\overline{w}=0$ we could call the resulting $\Gamma$ a curve of {\it zero width}, whatever that might mean. 

But in mathematics one should always trust the equations. It turns out that curves of zero width  have a long history as well as a  meaning. 

The width of a closed curve is defined as the distance between parallel tangents; it depends on the points of tangency.  Tangent lines are said to be parallel if the un-oriented lines they determine are parallel. It is possible for parallel but differently oriented tangents to determine the same un-oriented line, and in this case the distance between them is zero. A closed curve all of whose parallel tangents are of this kind can reasonably be said to have zero width.

Instantaneous width is non-negative. If the average width is 0, then the maximum instantaneous width is also 0, and therefore $w(t)=0$ for all $t$. So we have a Fourier series that represents zero. Therefore its coefficients are all zero, that is,
$$a_{2k}=b_{2k}=0 \quad \mbox{for all} \, k \geq 0$$
Hence the support function is a sum of odd harmonics, like all other cases of a curve of constant width: 
\be
\label{zero width}
p(t) =  \sum_{k=1}^{\infty} a_{2k+1} \cos (2k+1) t + b_{2k+1} \sin (2k+1) t \ee

This situation occurs when a simple curve, viewed from its interior, is everywhere concave except at cusps $P(t_0)$ where 
$$\lim_{t  \rightarrow (t_0)^{+}} T_{P(t_0)} = -\lim_{t \rightarrow (t_0)^{-}} T_{P(t_0)}$$
so the tangent vector $T_P$ instantaneously reverses direction, and the two limiting tangents lines are identical except for orientation, and there are no other parallel tangent lines. It is also possible when the curve is self-intersecting, in which case smooth concave arcs\footnote{That is, arcs whose curvature does not change sign.} join a collection of vertices where they still fail to have any parallel tangents: as a point progresses along the curve, the tangent turns in one direction. 

In the first publication to mention curves of constant width Euler  \cite{euler} considered a family of curves of zero width that he called {\it triangular curves} and showed how they could be used to construct curves of positive constant width, which he called  {\it orbiforms}. He then used the curves of constant width to construct catoptrices.

The triangular curves had three concave sides and three cusps. He called them {\it deltoids}; they probably reminded him of curved versions of the upper case Greek letter $\Delta$ (cp.  figure  \,\ref{Rab3_ZeroWidth&Evolute} which shows two such curves).

From our perspective,  a support function $p(t)$ with $w=0$, i.e. whose constant term is zero, represents the most general  curve of zero constant width. Since only odd harmonics appear, there are an odd number of cusps. If the series has more than one term, the curve may have self-intersections. 

Euler's method for constructing a curve of constant width from a triangular curve $\Delta$ is ingenious.\footnote{The description in White's otherwise valuable article \cite{white} is misleading.}  Although he presented his  construction for triangular curves there is nothing  that limits it to this case.
 Euler considers a line segment $\Lambda$ of fixed length $w$ at least equal in length to the longest arc of $\Delta$ ($w$ will turn out to be the width of the curve of constant width) tangent to $\Delta$.
 
  As the point of contact moves along $\Delta$ --- as the tangent line segment $\Lambda$ rolls along the curve without slipping ---  its endpoints  trace out a new curve. It is easy to see that each endpoint traces out three disjoint arcs that constitute one half of the curve, and that the endpoints of the arcs match up.  $\Lambda$ is perpendicular to the tangents  because  its endpoints have traced out the curve; hence it is always a diameter. Since it has length $w$, the curve has constant width.  Euler's process is related to the  geometrical notion of an {\it involute}, also {\it evolvent} (German: {\it evolvente}), which modern usage interprets as the locus of a taut thread unrolled against a curve, but here we see that the line segment --- the `taut thread' --- extends in both directions from the point of tangency, and that the `thread' --- better, the tangent line segment --- must have a fixed length.
  
Suppose $\Gamma$ is a curve of constant width $w$  with support function $p(t)$. The Fourier series  $p(t)-w/2$ is a support function that defines a curve of zero  width, say $\Delta$. Application of Euler's construction to $\Delta$ with a line segment of length $w$ reproduces $\Gamma$.

\subsection{Packing of curves of constant width}
\label{plane packing}

The rectangular packing of circles in the plane -- think of the projection of beer cans in a 24-pack carton -- has a density equal to  $\pi/4 \simeq 0.78539$: it covers only that fraction of the plane. Almost one-quarter of the space remains uncovered. Circles of fixed radius can be most efficiently packed in the plane if they are arranged on a lattice of regular hexagons. In 1773 Lagrange proved that the density of this cover,  $\pi/(2\sqrt{3} )\simeq 0.9060$, is maximal.

Consider the packing of Reuleaux curved triangles. Given a regular hexagonal lattice, three Reuleaux triangles can be arranged so that their vertices meet at the center of the hexagon. At this point the sum of the angles of the curved triangles is $2\pi$. The arrangement can be repeated indefinitely. This packing of the curved triangles covers $\pi/(3 \sqrt{3})+1/6 \simeq 0.77126$ of the plane. This is not as efficient as the rectangular packing of circles. 

However, just as the rectangular packing of circles is not optimal, this is not the optimal packing of Reuleaux curved triangles. There is a packing whose density  exceeds the density of the best packing of circles.  Perhaps cylinders with a Reuleaux base would be a preferable shape for a beer or soda can.

\thm There is a packing of Reuleaux curved triangles whose density is 
\be
\label{Reuleaux packing}
\frac{2 (\pi - \sqrt{3})}{\sqrt{15} + \sqrt{7} - 2 \sqrt{3}} \simeq 0.92288
\ee

\pf  It will suffice to cosnider the packing for Reuleaux curved triangles of width 2. 

The packing, which will be denoted $\mathfrak{R}_3$,  can be constructed as follows. Consider the rectangular lattice  $\mathfrak{L}=\big\{  (m (\sqrt{15}+\sqrt{7}-2\sqrt{3}),n ) : m,n \in \Z \big\}$. Choose as fundamental domain $\mathfrak{F}=\big\{ (x,y): 0 \leq x< \sqrt{15}+\sqrt{7}-2\sqrt{3}, -1  \leq y <1 \big\}$ and place a Reuleaux curved triangle of width 2 so that the midpoint of one of the circular arcs of  is tangent to the $y$-axis at $y=s/2$. Then place two congruent Reuleaux curved triangles so that the  midpoint of one the circular arcs is tangent to the line parallel to the $y$-axis  at $x= (\sqrt{15}+\sqrt{7}-2\sqrt{3})$, one at $y=-1$, the other at $y=1$, and so that each curved triangle intersects the interior of $\mathfrak{F}$. It will be seen that each of the parts of these two curved triangles is one-half a Reuleaux triangle and their union   is congruent to a Reuleaux triangle. Extend this pattern to the entire plane by translation by elements of the lattice. Figure \ref{R3 F domain} shows the fundamental domain $\mathfrak{F}$ and the parts of the packing of the Reuleaux curves that lie in it, and  figure  \ref{R3 packing} shows a section of the packing.

The density  in the plane of a Reuleaux triangle packing is equal to the ratio of the area they cover in $\mathfrak{F}$ to the area of $\mathfrak{F}$.  According to eq(\ref{area Reuleaux triangle}), the area of a Reuleaux triangle of width $2$  is $2(\pi - \sqrt{3})$.  Simple geometry shows that the  area of the rectangle $\mathfrak{F}$ is $2(\sqrt{15}+\sqrt{7}- 2\sqrt{3})$; thus the density is 
$$
\delta(\mathfrak{R}_3)=\frac{2 (\pi - \sqrt{3})}{\sqrt{15} + \sqrt{7} - 2 \sqrt{3}} \simeq 0.92288
$$
$\Box$

I conjecture that the packing by Reuleaux triangles has the greatest density for packings of the plane by a curve of constant width. Moreover, it seems likely that, given two curves of the same constant width, the one with the smaller area has the greater maximal packing density. The reason is that the curve that bounds the smaller area will be `less round' and offer more opportunities to snuggle closer to its neighbors. 

Another way to think about this is that the large symmetry group of the circle limits its ability to conform to its neighbors in a dense packing. Since the $n$-sphere has a large symmetry group, we are tempted to conjecture that the same will hold true in higher dimensions:  of two bodies in higher dimension of the same constant width, the one with the lesser volume will admit a packing with the greater density. Thus, in 3-dimensions, if, as is generally expected, it turns out that the two Meissner bodies minimize the volume of solids of given constant width, then our conjecture implies  they will have the same maximal packing density (because they have the same volume) and it will be  the greatest possible packing density.

\section{Surfaces of constant width}
 
\subsection{Initial considerations}
 
Few solids of constant width are explicitly known.\footnote{In analogy to Euler's {\it orbiforms}, solids of constant width are sometimes called {\it sphaeriforms} \cite{bonnesen+fenchel} or {\it spheriforms}.}  Can similar methods be applied to construct  them? One might believe that Reuleaux' construction of the curved polygons of constant width could be generalized to the regular solids. Consider, for instance, the solid defined by the intersection of four  balls of radius $w$ with centers at the vertices of a regular tetrahedron of edge $w$. Each point of its surface $\Sigma$  is at distance $w$ from the opposite vertex. Is this not a surface of constant width? It is not. The distance between the midpoints of opposite edge arcs is $\big( \sqrt{3}-\frac{\sqrt{2}}{2} \big) \, w \sim 1.0249 w$.
 
 Meissner \cite{meissner1, meissner2, meissner+schilling} showed that the Reuleaux tetrahedron could be modified to make it of constant width by `rounding' three of the edge arcs. There are two ways to choose 3 edges of a tetrahedron -- sharing a vertex, or not -- so there are two Meissner tetrahedra.

 If $\Gamma$ is a planar curve of constant width with an axis of symmetry  $AA'$, then the surface of revolution obtained by rotating $\Gamma$ about the axis $AA'$ is a solid of constant width. This is almost self-evident.  The normal to the surface at any point is directed toward the axis of revolution and the intersection of the plane determined by these lines with the solid is a curve congruent to $\Gamma$.  Figure \ref{HLR_S0} shows perhaps the simplest example, corresponding to the support function $p(t,u)=1+\frac{1}{8} \cos 3t$. The cross sections through the axis of symmetry are the curve of constant width shown in  figure   \ref{Rab3}. Figure \ref{HLR_S0_hue} shows two copies the same surface, color coded to display the the two principal radii of curvature. On the left we see $ \rho_0(t,u)$. The radius of curvature at the vertex is 0; at the base, 2. In between, the cyan band corresponds to radius 1. The radius of curvature $ \rho_1(t,u)$ is shown on the right.
 
Apart from this general construction and Meissner's tetrahedra, I have been unable to find explicit examples of solids of constant width. This of course provides motivation to construct new examples.


\subsection{Ovoid surfaces}

An oval surface is usually called an \index{ovoid}{\it ovoid}. Like the oval, there is no accepted definition of an ovoid. We shall say that an ovoid $\Sigma$ is a bounded simply connected convex smooth surface that has a tangent plane at each point, and that $\Sigma$ lies entirely on one side of the tangent plane.\footnote{David Pollen brought to my attention that an unbounded tube over a curve of constant width is an unbounded surface of constant width. The compactness condition excludes such surfaces.}A surface of constant width is a special case of an ovoid.



A simply connected smooth surface is convex fif the curve $\Gamma$ that is the intersection of the ovoid $\Sigma$ with a plane through the normal at $P$ has non-negative curvature for every point $P \in \Sigma$  and plane through the normal at $P$. According to a famous theorem of Euler, the curvature $\kappa$ of such a curve $\Gamma$ is a linear combination of the two principal curvatures at $P$ -- the maximum and minimum curvatures in the various directions at $P$ -- and has the form $\kappa_{\theta} = \kappa_0 \cos^2  \theta+ \kappa_1 \sin^2 \theta$ where $\theta$ denotes the angle between the direction of interest and the direction of the principal curvature $\kappa_0$. Hence $\Sigma$ is convex fif the minimum principal curvature is non-negative at every point.

The principal radii of curvature at $P$ in the directions of the principal curvatures are defined by
\be
\rho_k = 1/\kappa_k, \quad k=0,1
\ee

\subsection{Constant width condition}
Consider a cartesian coordinate system and a point $P =(x,y,z)$ on a hypothetical ovoid surface $\Sigma$. This means that the principal radii of curvature are non-negative at every point of $\Sigma$. We seek conditions that the ovoid will be a surface of constant width $w$. Assume the origin $O$  lies inside $\Sigma$ and that $p$ denotes the distance from the origin to the plane tangent to $\Sigma$ at $P$. Following the approach taken in section \ref{support p}, introduce spherical coordinates for the support function $p$. The unit normal to the tangent plane is the vector

$$\mathbf{n} = (\cos t \sin  u, \sin t \sin u, \cos u) $$
where $t, \, 0 \leq t < 2\pi $, denotes the angle measured from the positive $x$-axis and $u, \, 0 \leq u < \pi$, measures the angle from the positive $z$-axis.  

Parametrize the coordinates of $P$ by the angles $(t,u)$: $x=x(t,u)$ etc.; note that  these are not the spherical coordinates of $P$. Then

\be
\label{p=xyz}
p(t,u) = x\, \cos t \sin u  + y\, \sin t \sin u + z\, \cos u
\ee
and the direction opposite to the normal $\mathbf{n}$ is given by the expression
$$t \rightarrow t+\pi, \quad u \rightarrow \pi - u $$
Hence the instantaneous width of $\Sigma$ along the line determined by  $\mathbf{n}$ is
\be
w(t,u) = p(t,u) + p(t+\pi, \pi-u) 
\label{cw identity}
\ee

Now consider the general Fourier series representation of $p(t,u)$:
\bea
p(t,u)&=& \sum_{m,n \geq 0} \big(
 cc_{mn} \cos m t \cos n u + 
 sc_{mn} \sin m t \cos n u  + \nonumber  \\
 && \qquad \qquad  
 cs_{mn} \cos m t \sin n u  +  
 ss_{mn} \sin m t \sin n u   \big)
\eea
and suppose that the instantaneous width is constant: $w(t,u)=w$. Then
\bea
w &=& p(t,u)+ p(t+\pi, \pi-u)  \nonumber  \\
 &=& \sum_{m,n \geq 0} \Big( 
 \big(1 + (-)^{m+n} \big)  \big( cc_{mn} \cos m t \cos n u + sc_{mn} \sin m t \cos n u \big) +   \nonumber \\
 && \qquad  \big(1 - (-)^{m+n}\big) \big( cs_{mn} \cos m t \sin n u  +    ss_{mn} \sin m t \sin n u   
 \Big)
\eea
The right side reduces to the constant term, which yields
\be
cc_{00} = \frac{w}{2}
\ee
Suppose $m+n>0$. Then  $cc_{mn} = sc_{mn} =0$ if $m+n$ is even and $cs_{mn}=ss_{mn}=0$ if $m+n$ is odd. Thus candidate width functions have the form
\bea
\label{p solid}
p(t,u) 
 &=& c_{00}+ \nonumber \\
 &&  \frac{1}{2 } \sum_{n> 0}  \sum_{m \geq 0} \Big\{
 \big(1 - (-)^{m+n} \big)  \big( cc_{m,n-m} \cos m t \cos (n-m) u + \nonumber \\
 && \qquad \qquad sc_{m,n-m} \sin m t \cos  (n-m) u \big) +   \nonumber \\
 && \qquad  \big(1 + (-)^{m+n}\big) \big( cs_{m,n-m} \cos m t \sin  (n-m) u  +  \nonumber \\
 && \qquad \qquad  ss_{m,n-m} \sin m t \sin  (n-m) u    \big)
 \Big\}
\eea
with $c_{00}=w/2$. Any convergent Fourier series that satisfies these conditions satisfies the constant width identity eq(\ref{cw identity}). 

Equations for $x,y,z$ can be found by requiring that $p$ be an extremum at $(x,y,z)$, that is, 
$$p_t = p_u = 0 $$
where $p_t = \partial p/ \partial t$ and $p_u = \partial p/ \partial u $ which, when applied to eq(\ref{p=xyz}),  leads to
\bea
x &=& \big(  p_u   \cos u   + p \sin u \big) \cos t  -  p_t   \frac{ \sin t }{\sin u}  \nonumber \\
 y &=& \big(  p_u   \cos u   + p \sin u \big)  \sin  t +  p_t   \frac{ \cos t }{\sin u}  \nonumber \\
 z &=&  p \, \cos u - p_u \sin u 
 \eea
 
 This implies, in particular, that the square of the distance from $OP$ is
 $$ p^2 + p_u^2 + \Big(\frac{p_t}{\sin u} \Big)^2, $$
 which reminds us of the problems inherent in using spherical coordinates: they  are not invertible and can lead to singularities at the poles where $u=0$ or $u=\pi$. Since the points of a surface of constant width are at finite distance from $O$ it follows that the derivative $p_t$ must have some factor whose behavior is proportional to that of  $\sin u$ where the latter has zeros. 
 
Assume, for instance, 
 \be
\frac{\partial  p(t,u)}{\partial t} =q(u,t)  \sin^k u 
  \ee
 for some positive integer $k$. This enforces finiteness of the distance $PO$ if the derivatives of $p$ are finite.\footnote{It will turn out that $k \geq 2$ will be necessary. 
The factor $\sin u$ in $\partial p(t,u)/\partial t$ cancels the same factor in the denominator of the coordinates $(x(t,u),y(t,u))$  but the result, although finite, may not have a well-defined limit as $u \rightarrow 0$ or $\pi$. The circles of latitude are squeezed to a point for these limiting values. An additional factor of $\sin u$ in $\partial p(t,u)/\partial t$  insures that $p(t,0)$ and $p(t,\pi)$ will be well-defined.}
 
 Suppose that
 $$ p(t,u) = \frac{w}{2} + f(u) +   g(t,u) \sin^k u  $$
The width becomes, from eq(\ref{cw identity}),
 \be
 \label{w2}
 0 = \big(  f(u)+f(\pi-u) \big) + \big(  g(t,u)+g(t+\pi,\pi-u) \big) \sin^k u
\ee
  whence 
 \be
\label{g}
g(t,u) + g(t+\pi, \pi -u) =0
\ee
and
\be
\label{f}
  f(t,u)+f(t+\pi,\pi-u)  =0,
\ee
which is just the condition that the curve defined by the intersection of $\Sigma$ with the plane $t={\rm constant}$ is a curve of constant width $w$. Notice that the form of the fundamental equation  \ref{p solid} for the support function is maintained regardless of the parity of $k$.

These constraints on the form of $p(t,u)$ are necessary but not sufficient that $\Sigma$ be a surface of constant width because it may not be convex. Convexity can be enforced by insuring that the principal curvatures (equivalently, the principal radii of curvature) be positive for all parameter pairs $(t,u)$. The principal radii can be calculated from the first and second fundamental forms of $\Sigma$ in the usual way. The expressions quickly become analytically unmanageable except for computer calculation, although certain general observations can be made. Here is one of them.

Let $\Sigma$ be a surface of constant width $w$. Consider points $P,Q \in \Sigma$ and that $P=(x(t,u), y(t,u), z(t,u))$. We shall say that $P$ and $Q$ are opposite if $Q=(x(t+\pi,\pi-u),y(t+\pi,\pi-u),z(t+\pi,\pi-u))$. Then the tangent planes at $P$ and $Q$ are parallel and the distance between them is $W$. Write $\rho_i(P)$ for the $i$-th principal curvature at $P$.

\thm  \label{surface curvature thm} Suppose that  $\Sigma$ is a surface of positive constant width $w$ and $P,Q \in \Sigma$ are opposite points. Then
\be
 \rho_0(P) +  \rho_1(Q) = w
\ee

\pf   The expressions for $ \rho_0$ and $ \rho_1$ depend on  the support function $p$ and its derivatives through second order at $P$ and $Q$.  The proof of the first equation consists of simplifying it, taking into account that $p(t,u)+p(t+\pi,\pi-u)=w$.
 \begin{center}
\rule{60pt}{0.5pt}
  \end{center}
  
\begin{small}
Some readers may be interested in the details.  The radii of curvatures are roots of a quadratic polynomial, so we may write
$$\rho_0 = A + \sqrt{B}, \quad \rho_1 = A-\sqrt{B} $$
where $A$ and $B$ are explicitly given in terms of the support function by
\bea
A(t,u) &=& \frac{1}{2} \left(p^{(0,2)}(t,u)+\cot (u)
   p^{(0,1)}(t,u)+\csc ^2(u) p^{(2,0)}(t,u)+2
   p(t,u)\right) \nonumber \\
   && \label{A}\\
B(t,u) &=&\frac{1}{16} \csc ^4(u) \Big( \sin ^2(2 u)
   (p^{(0,1)})^2+  \nonumber \\
   && 8 \sin (u) \cos (u)
   p^{(0,1)} \big(p^{(2,0)}- \sin ^2(u)
   p^{(0,2)}\big)+  \nonumber \\
   &&4  \big(  (p^{(2,0)})^2+\sin ^4(u)
   (p^{(0,2)})^2-2 \sin ^2(u) p^{(0,2)}
   p^{(2,0)}+  \nonumber \\
   && 4 \big(\cos (u)
   p^{(1,0)}-\sin (u)
   p^{(1,1)}\big)^2\big) \Big)
   \label{B}
\eea
 where $p^{(m,n)}= \partial^{m+n}p(t,u)/\partial^mt \partial^n u$.

In order to simplify 
$$
\rho_0(t,u) +\rho_1(\pi+t,\pi-u) =  A(t,u)+A(\pi+t,\pi-u) +\sqrt{B(t,u)}-\sqrt{B(\pi+t,\pi-u)}
$$
we need expressions for the derivatives at $(\pi +t,\pi-u)$. From $p(t,u)+p((\pi +t,\pi-u)=w$ find
\bean
p(\pi+t,\pi-u) &=& w - p(t,u) \\
p^{1,0}(\pi+t,\pi-u) &=& -p^{1,0}(t,u) \\
p^{1,1}(\pi+t,\pi-u) &=&  p^{1,1}(t,u) \\
p^{2,0}(\pi+t,\pi-u) &=& -p^{2,0}(t,u) \\
p^{0,1}(\pi+t,\pi-u) &=& p^{0,1}(t,u) \\
p^{0,2}(\pi+t,\pi-u) &=& -p^{0,2}(t,u) \\
\eean
Substitution of these relations in $A(\pi+t,\pi-u)$ and $B(\pi+t,\pi-u)$  yields
\bean
A(t,u)+A(\pi+t,\pi-u) &=& p(t,u)+p(\pi+t,\pi-u) =w \\
 B(t,u) &=& B(\pi+t,\pi-u)
 \eean
which completes the proof.
\end{small}\\
$\Box$



By interchanging $P$ and $Q$ and adding the results, the theorem implies
\be
\label{mean radii condition}
\rho_{{\rm mean}}(P)  +\rho_{{\rm mean}}(Q)  = w
\ee
This equation simply states that the sum of the mean curvature at $P$ and $Q$ is the width. Although the expression of each principal radius in terms of the support function involves a quadratic irrationality --- a square root of a complicated expression in the derivatives of $p$ of order $
\leq 2$, it turns out that the corresponding expression for the sum of the mean radii at opposite points in eq(\ref{mean radii condition}) is much simpler. In fact, 
$$ \sin^2u \left( \rho_{{\rm mean}}(P) +   \rho_{{\rm mean}}(Q)  - w \right)$$
is again a Fourier series whose coefficients are linear combinations of the original coefficients of $p(t,u)$. This observation leads to the following

\thm  
\label{thm ovoid} Suppose that $\Sigma$ is an ovoid with support function $p$ and that $P$ and $Q$ are opposite points and that $w$ is a positive constant. Then  $\Sigma$ is a surface of constant width $w$ fif 
$$ \rho_{{\rm mean}}(P) +   \rho_{{\rm mean}}(Q)  = w$$

\pf According to eq(\ref{mean radii condition}),  if $\Sigma$  is a surface of constant width $w$, then the mean radii of curvature satisfy the condition of the theorem.

We will outline the proof of the converse, suppose that\footnote{Note that $ \sin^2(\pi-u)=\sin^2u$.}
\be
\label{thm mean radii}
\sin^2u \left( \rho_{{\rm mean}}(P) +   \rho_{{\rm mean}}(Q)  - w \right)=0
\ee
Since the convexity condition is incorporated in the hypothesis, it suffices to check that the Fourier series for $p(t,u)$ has the right form, that is, the form eq(\ref{p solid}).  This will follow from consideration of the above Fourier series  
whose coefficients must all vanish. What are those coefficients?  Because the expression is linear in the coefficients of  $p$ it will be enough to consider the components independently. Consider the summands of $p$ that have the form
\bean
cc_{m,n-m} \cos(n-m)u \cos m u + cs_{m,n-m} \cos(n-m)u \sin m u + \\
sc_{m,n-m} \sin(n-m)u \cos m u + ss_{m,n-m} \sin(n-m)u \sin m u 
\eean
where $0 \leq m   \leq n$ for non-negative $n$. These summands for the mean curvature include various harmonics in $u$ less than and greater than $(k-m)u$ that overlap with the expressions for adjacent values of $k$ and $m$, but the  terms involving $(k-m)u$ appear only once as coefficients of the orthogonal functions $\cos (m t \pm (k-m)u)$ and $\sin (m t \pm (k-m)u)$. Their coefficients yield the equations 
\begin{small}
\bean
0&=&-\frac{1}{8} \left(k^2-2 k m+3 m^2-2\right)
   \left(\left((-1)^k+1\right)
   \text{cc}(m,k-m)+\left((-1)^k-1\right)
   \text{ss}(m,k-m)\right) \\
0  &=&  -\frac{1}{8} \left(k^2-2 k m+3 m^2-2\right)
   \left(\left((-1)^k+1\right)
   \text{cc}(m,k-m)-\left((-1)^k-1\right)
   \text{ss}(m,k-m)\right)\\
 0  &=& \frac{1}{8} \left(k^2-2 k m+3 m^2-2\right)
   \left(\left((-1)^k-1\right)
   \text{cs}(m,k-m)-\left((-1)^k+1\right)
   \text{sc}(m,k-m)\right) \\
   0 &=& -\frac{1}{8} \left(k^2-2 k m+3 m^2-2\right)
   \left(\left((-1)^k-1\right)
   \text{cs}(m,k-m)+\left((-1)^k+1\right)
   \text{sc}(m,k-m)\right)
   \eean
\end{small}
Unless one or another of the factors $\big( 1 \pm (-1)^k \big) =0$, the solution of these equations is
$$
cc_{m,k-m}=0, \quad cs_{m,k-m}=0, \quad sc_{m,k-m}=0, \quad ss_{m,k-m}=0
$$
When the factors are zero, there is no constraint on the corresponding coefficients of the Fourier series for the support function. These are just the coefficients that must be zero in order for $p$ to be the support function of a surface of constant width.\\
   $\Box$

\subsection{Constant width surfaces of revolution}

The simplest way to satisfy this condition is to take $g(t,u)=0$; this leads to the family of surfaces of revolution of constant width obtained by rotating the curve of constant width defined by $f(u)$ about the $z$-axis. In this case, calculation of the principal radii (the reciprocals of the principle curvatures) yields:
\bea
\label{principal radii}
 \rho_0(t,u) &=& 1+ f(u)+f''(u) = p + p'' \nonumber \\
 \rho_1(t,u) &=& 1 + f(u) + f'(u) \,\cot u = p + p'\, \cot u
\eea

Of course $f$ must be constrained for the surface to be convex, just as was required for a curve a constant width. In this case both principal radii of curvature must be positive, that is, $f$ must be a function such that  $ \rho_0 \geq 0,\quad  \rho_1 \geq 0$. We recognize the first equation as  eq(\ref{p radius of curvature}) for the radius of curvature of a curve of constant width -- here, the profile curve of the surface of revolution. So this is just what one expects. The second equation appears to be supplemental. But here, once again, we see the consequences of the singularities of spherical coordinates. The variable $u$ is constrained to $0 \leq u < \pi$ and the constant width condition eq(\ref{f}) when applied to a Fourier series for $f$ eliminates all terms of the form $\sin n u$. This implies that $\cot u \sin nu$ is a polynomial in $\cos u$ and $\sin u$, hence bounded.  Judicious choice of the Fourier coefficients insures convergence and positivity, but these constraints are already subsumed in the positivity condition for the other principal radius.

\subsection{New examples}
\label{surface examples}
Other choices lead to entirely new surfaces of constant width. An interesting example is given by $f(u) = (1/16)\cos 3 u$, $g(t)= (1/16) \cos 3 t$, that is,
\be
\label{HLR_S1}
\Sigma: \quad p(t,u) = 1 +  \frac{1}{2} \Big( \frac{1}{8}\cos 3u +  \frac{1}{8}  \cos 3 t \,  \sin^2 u \Big)
\ee
The surface, illustrated in  figure  \,\ref{HLR_S1}, is `tetrahedral-like' but it is not Meissner's surface.

Generally speaking, the shape of $\Sigma$ is roughly determined by the lowest order frequencies in the Fourier series for $p(t,u)$. 

If $p(t,u)=w/2$ ,  $\Sigma$ is a sphere of diameter $w$. Other terms in the Fourier series modify the sphere while maintaining constant width.

Since $0 \leq t <2\pi$, if the frequencies are  odd numbers, say $2k+1$, the solid will have $2m+1$ `sides' ranged symmetrically about the  polar axis $u=0$. Since $0 \leq u < \pi$, the solid will have  $2n+\frac{1}{2}$ `sides' ranged between the pole $u=0$ and the 'pole' $u=\pi$, which the partial side touches.

If $t$ does not appear in the series, then $\Sigma$ is a surface of revolution. The profile is a symmetric curve of constant width $w$.   If only $\cos u$ or $\sin u$ appear, the surface is the sphere.

Summands that contain $t$ must be multiplied by $\sin^n u$ for some $n\geq 2$; this suppresses the singularity of the spherical coordinate system at $u=0$. Thus an expression of the form
$$ \Sigma : \quad  p(t,u) = 1 + a \cos (2m+1)t \, \sin^n u $$
will produce a surface which, if it is convex, will have $2m+1$ `sides' organized as lunes extending from one pole to the other. Convexity will be assured by selecting $a$ sufficiently small. 

The simplest case of this situation is
\be
\Sigma : \quad p(t,u) =  1 +\frac{1}{2} \cos t \, \sin^2 u
\ee
$\Sigma$ is not a surface of revolution (cp.  figure  \,\ref{HLR_S(1,0)}); this oddball surface of constant width has  but one `side' as $t$ varies.

The expression
\be
\label{HLR_S(5,3)}
\Sigma : \quad p(t,u) = 1 + \frac{1}{2} \Big( \frac{1}{8} \cos 3 u + \frac{1}{24} \cos 5 t \,  \sin^2 u  \Big)
\ee
produces a surface of constant width with pentagonal and triangular symmetries ( figure  \,\ref{HLR_S(5,3)}).

Although the condition of constant width is very constraining and causes the surfaces to look like deformations of the sphere (which they are), some unusual configurations can be built. Consider, for instance, 
$$ p(t,u)=1+ \frac{1}{10}\Big(  \cos 3t \sin 3u \sin^2 u\Big) $$
 The support function can also be written as
\bean
p(t,u) &=&1+ \frac{1}{80} \Big( \sin (3t+u) - \sin (3t-u) +2 \sin (3t+3u)  -  \\
&& \qquad   2 \sin (3t-3u) + \sin (3t+5)  - \sin (3t-5u) \Big) 
\eean
The surface,  shown in  figure   \ref{HLR_S(3,3)}, appears to be partitioned into 9 curved quadrilaterals. For $ \rho_0$ the minimum radius of curvature is 0 and the maximum is 2/10; for $ \rho_1$ the maximum is 2/10 and the maximum is 2.

\subsection{Area and Volume} Blaschke \cite{blaschke2} proved a remarkable formula  relating the area   $A_\Sigma$  of a surface of constant width $w$ and the 3-dimensional volume $V_{\Sigma}$ it bounds:
\be
\label{blaschke eq}
V_{\Sigma} = \frac{w}{2} A_{\Sigma} - \frac{\pi}{3} w^3
\ee
Contrast this with the theorem of Barbier for plane curves of constant width which may bound different areas although  all of them have the same perimeter.
 
The Meissner surfaces have the same area $ A_{M}$, hence the same volume $V_{M}$ \cite{kawohl+weber}: 
 \bean
 A_{M} &=&  \Big( 2 -\frac{\sqrt{3}}{2} \arccos \frac{1}{3}  \Big) \pi\, w^2 \simeq 2.934115 \, w^2 \\
V_{M} &=& \Big( \frac{2}{3} -\frac{\sqrt{3}}{4} \arccos \frac{1}{3} \Big) \pi\, w^3 \simeq 0.41980 \, w^3
\eean

Two years after Meissner published his construction, Bonnesen and Fench\-el \cite{bonnesen+fenchel} speculated that the Meissner solids realize  the smallest volume bounded by any surface of fixed constant width. This conjecture is still open although alternatives have been increasingly tightly circumscribed. 

We close this paper with the formula for calculating the volume given  the support function.  The volume element in spherical coordinates is 
$$ r^2 \sin \phi \, dr \, d \phi \, d \theta \quad \mbox{where} \quad 0 \leq 0 < \phi, 0 \leq \theta < 2 \pi  $$
but recall that the angles $(t,u)$ are not the angles in the spherical coordinate representation of $P \in \Sigma$. The connection is established by
\bean
r(t,u)&=&\big( x(t,u)^2 + y(t,u)^2 + z(t,u)^2 \big)^{1/2}\\
\phi(t,u) &=& \arccos z(t,u)/r(t,u) \\
\theta(t,u) &=& \arctan y(t,u)/x(t,u)
\eean
from which the volume element can be expressed in terms of $t$ and $u$ as (The subscripts denote differentiation)
$$ dr \wedge d\phi \wedge d \theta =\big ( \phi_u \theta_t - \phi_t \theta_u \big)  \, dr \wedge du \wedge d t  $$
Therefore the volume of the solid bounded by $\Sigma$ is
\bea
V_{\Sigma} &=& \int_0^{2\pi} \int_0^{\pi} \int_0^r r^2 \sin \phi  \, dr\,   d\phi  \, d \theta \nonumber  \\
&=&\frac{1}{3}   \int_0^{2\pi} \int_0^{\pi}  r^3 \,   \big ( \phi_u \theta_t - \phi_t \theta_u \big)  \sin \phi  \,  du  \,d t
\eea

\section{Comments and conclusions}

\begin{itemize}

\item It would be nice to find a formulation for surfaces that avoids the problems of coordinate system singularities while retaining the calculational advantages of Fourier series.

\item In all dimensions, the question of packing bodies of constant width remains to be studied. The result for the Reuleaux curved triangle given here is but a small step in this direction, since it has not been proved that this packing has the greatest possible density.
\end{itemize}

What is new in this paper?  Insofar as I know,
\begin{enumerate}
\item The  expression for the curvature of an oval eq(\ref{p convexity}).
\item The inequality for the radius of curvature of an oval eq(\ref {oval m rho bound}).
\item The generalization of Barbier's theorem to ovals (section \ref{oval Barbier thm}).
\item For a curve of constant width, that the sum of the maximum and minimum radii of curvature is equal to the width (section \ref{rhoMin+rhoMax}).
\item  The explicit  density for a packing of Reuleaux curved triangles that exceeds the density for circles give by eq(\ref{Reuleaux packing}).

\item The conjectures about packing density for curves and bodies of constant width that follow  eq(\ref{Reuleaux packing}).

\item  For surfaces of constant width $ \rho_0(P)+ \rho_1(Q)=w$ where  the $\rho_i$ are the  principal curvatures, $P$ and $Q$ are opposite points, and  $w$ is the width (eq(\ref{surface curvature thm})).

\item An ovoid is a surface of constant width fif the sum of the mean radii of curvature at opposite points is equal to the width (section \ref{thm ovoid}).

\item And, of course, the variety of new examples of curves and surfaces of constant width.

\end{enumerate}

\newpage

\section{Figures}


\begin{figure}[h]
\begin{center}
\includegraphics[width=2.5in]{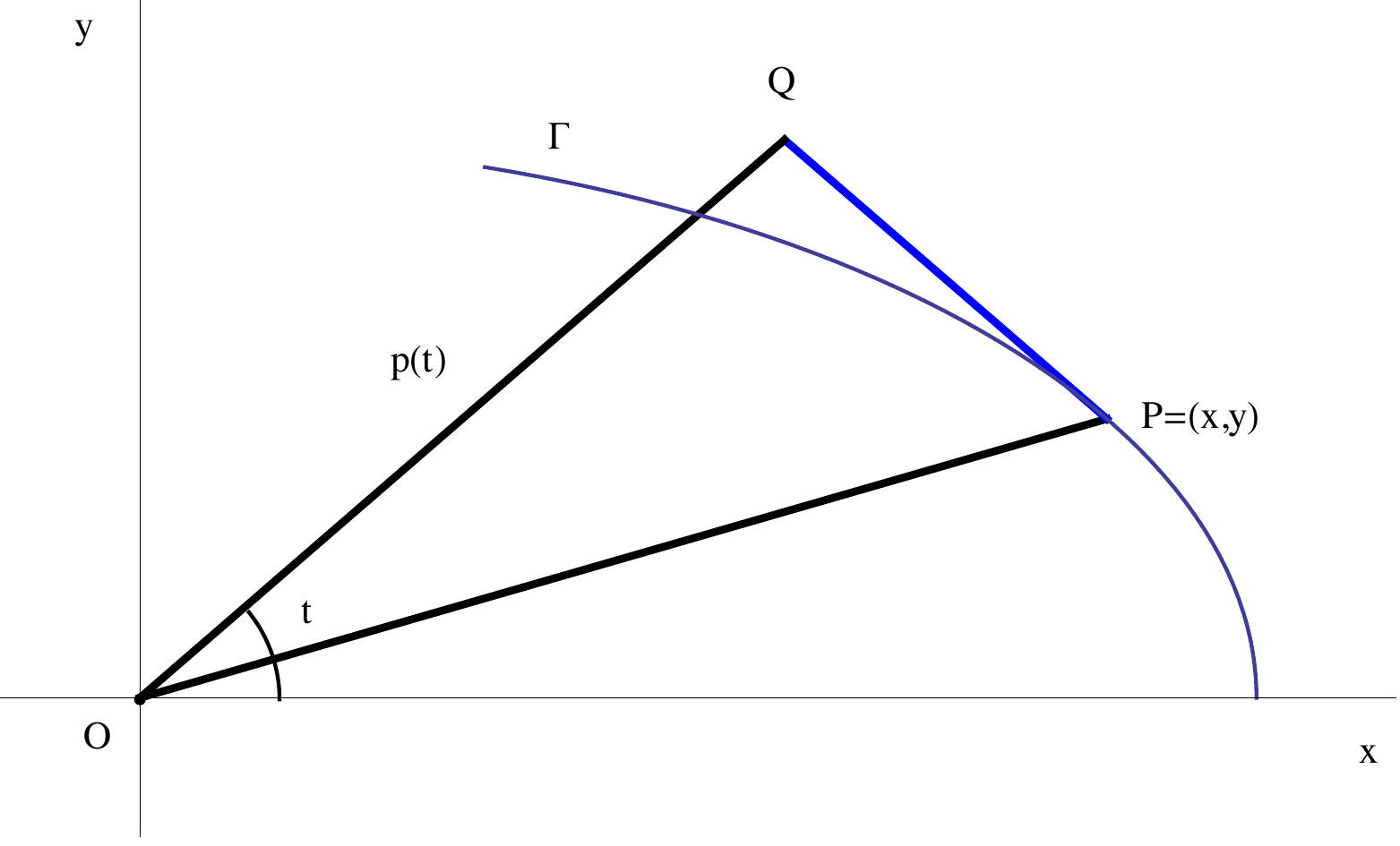}
\caption{Geometry of support function $p(t)$.}
\label{geometry}
\end{center}
\end{figure}

\begin{figure}[h]
\begin{center}
\includegraphics[width=1.25in]{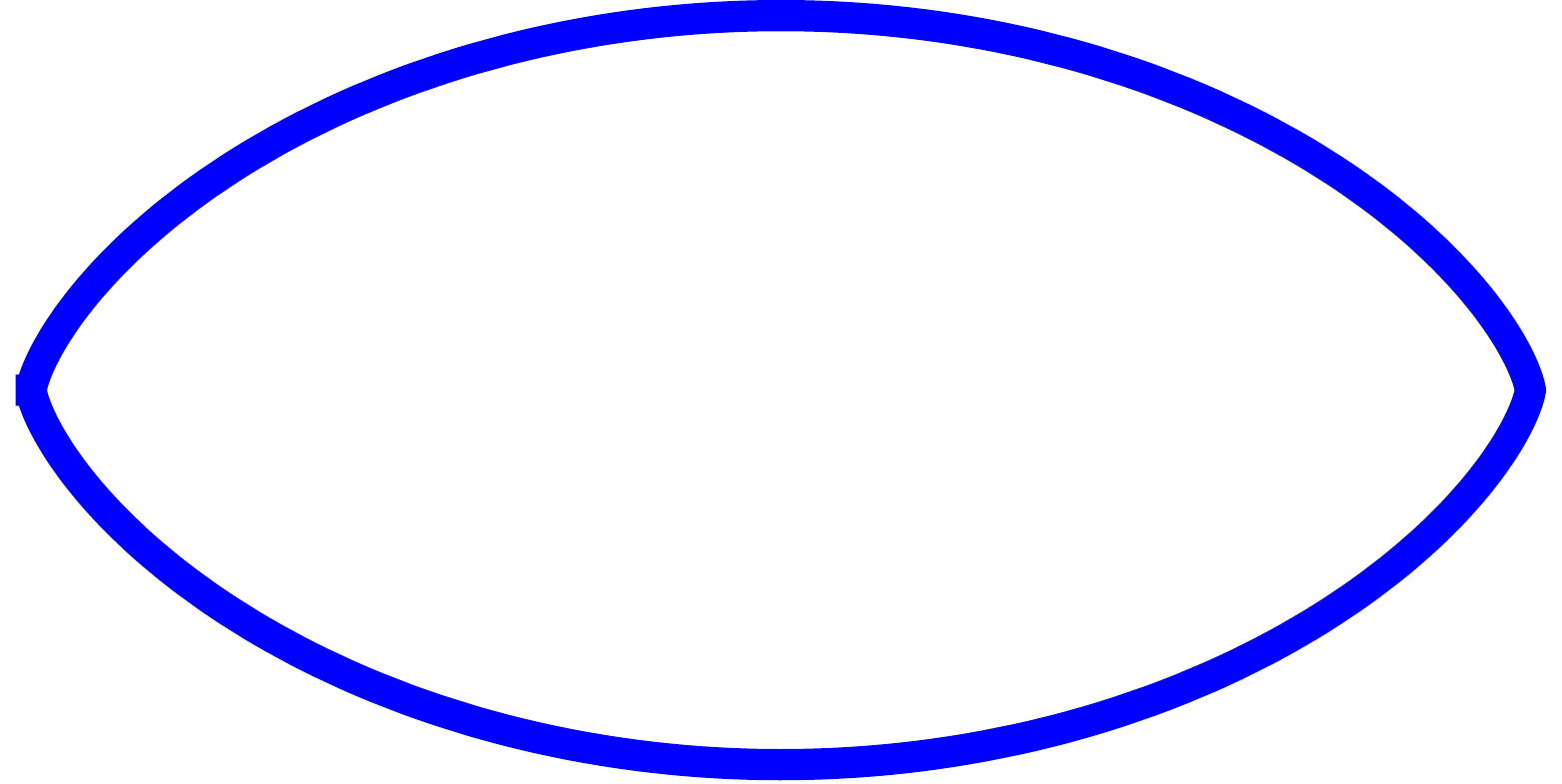}
\caption{Lozenge shaped oval of degree 2 with support function $p(t)=1 +\frac{1}{3}\cos 2 t$.}
\label{fig: deg 2 lozenge}
\end{center}
\end{figure}

\begin{figure}[h]
\begin{center}
\includegraphics[width=1.in]{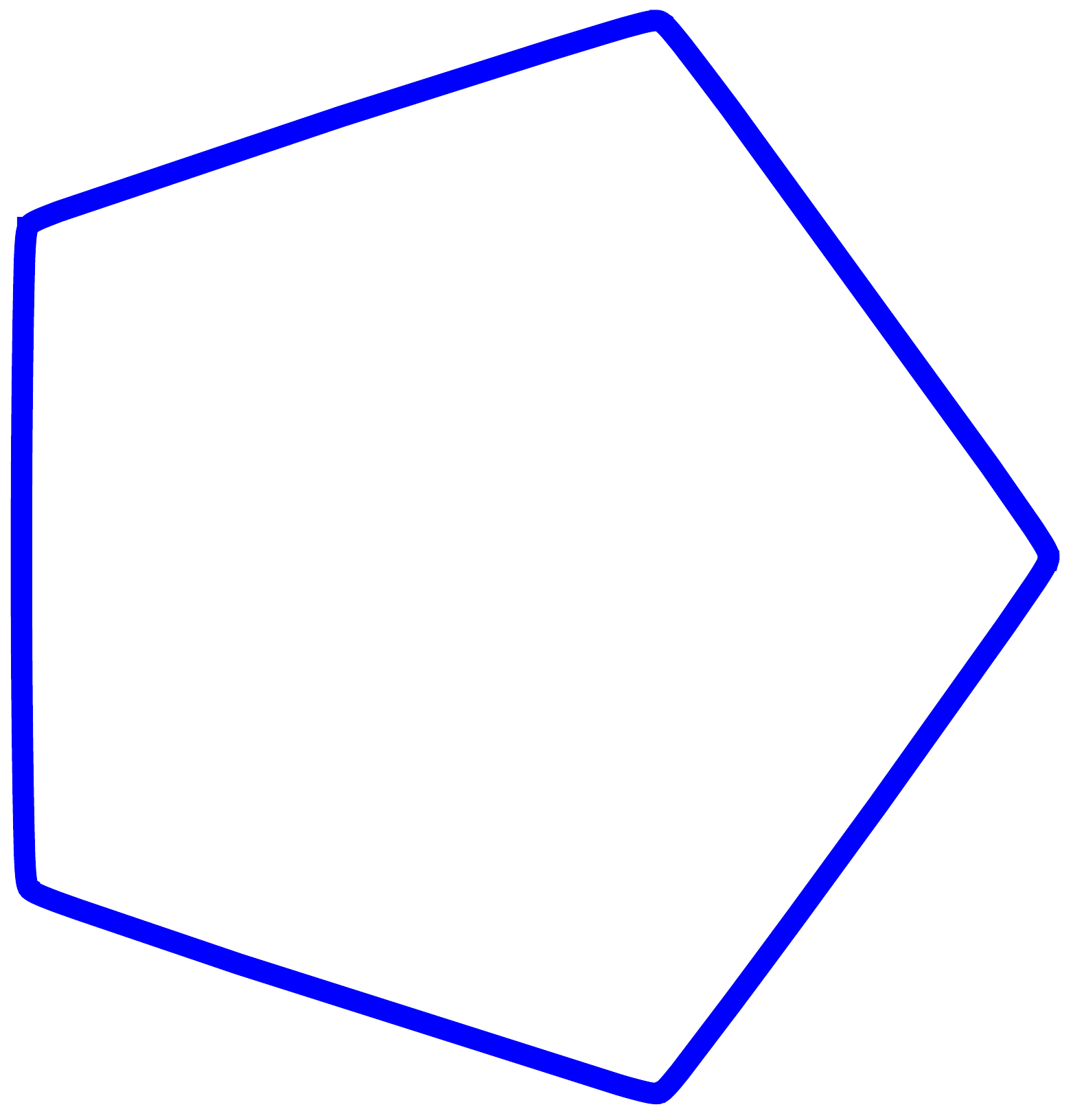}
\caption{Oval with maximum radius of curvature 32 and pentagonal symmetry.}
\label{fig: Fejer n=32,sigma=5}
\end{center}
\end{figure}


\begin{figure}[h]
\begin{center}
\includegraphics[width=1.5in]{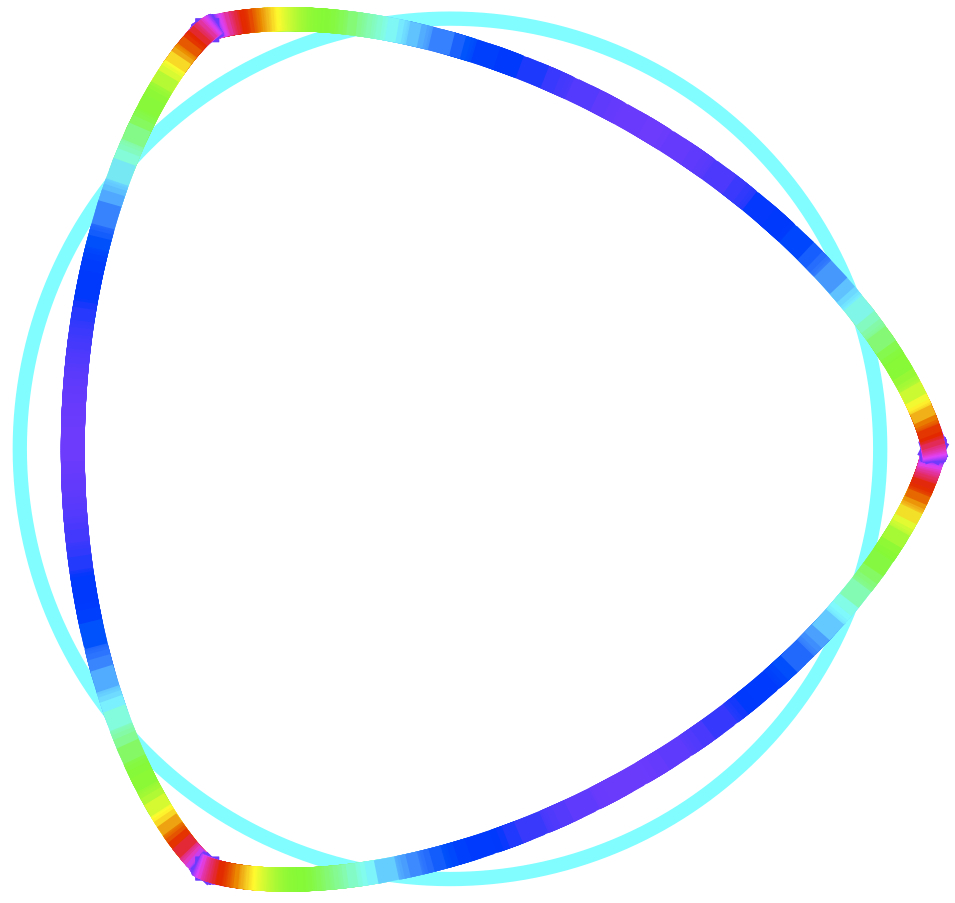}
\caption{Constant width curve for  $p(t)=1+(1/8) \cos 3t$  and circle of same width. Color codes the radius of curvature: $0<{\rm red}<{\rm cyan}=1 < {\rm blue}<	2$.}
\label{Rab3}
\end{center}
\end{figure}

\begin{figure}[h]
\begin{center}
\includegraphics[width=1.5in]{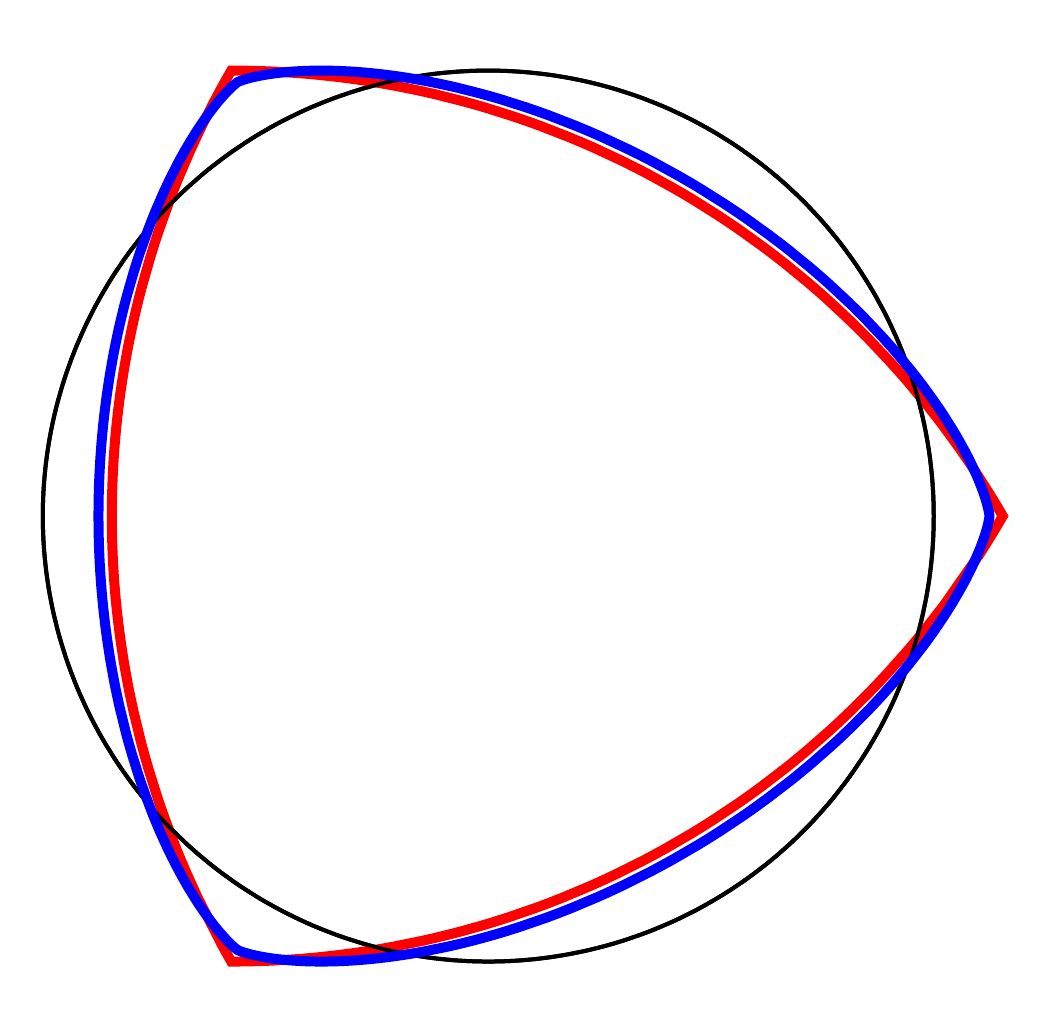}
\caption{Reuleaux's curved triangle (red) and curve of same width for $p(t)=1+\frac{1}{8}\cos 3 t$ (blue).}
\label{HLR3+reu3}
\end{center}
\end{figure}

\begin{figure}[h]
\begin{center}
\includegraphics[width=1.5in]{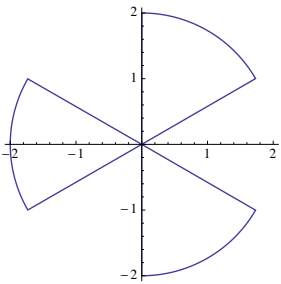}
\caption{Polar plot of radius of curvature for Reuleaux's curved triangle as function of parameter $t$.}
\label{Reu3_rho(t)_polarplot}
\end{center}
\end{figure}

\begin{figure}[t]
\begin{center}
\includegraphics[width=1.5in]{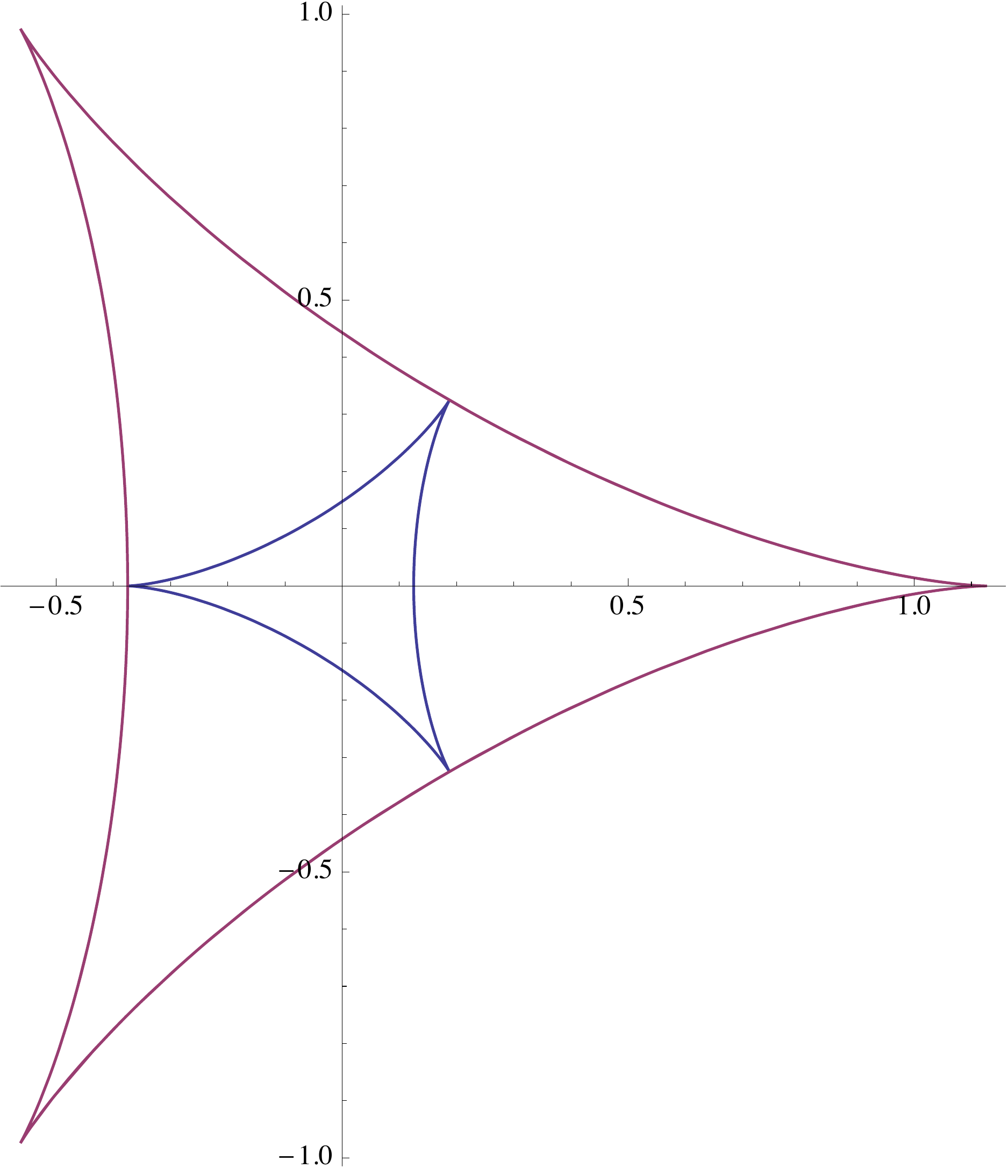}
\caption{Curve of zero width (inner) and its evolute for $p(t)=(1/8) \cos 3t$.}
\label{Rab3_ZeroWidth&Evolute}
\end{center}
\end{figure}

\begin{figure}[h]
\begin{center}
\includegraphics[width=1.5in]{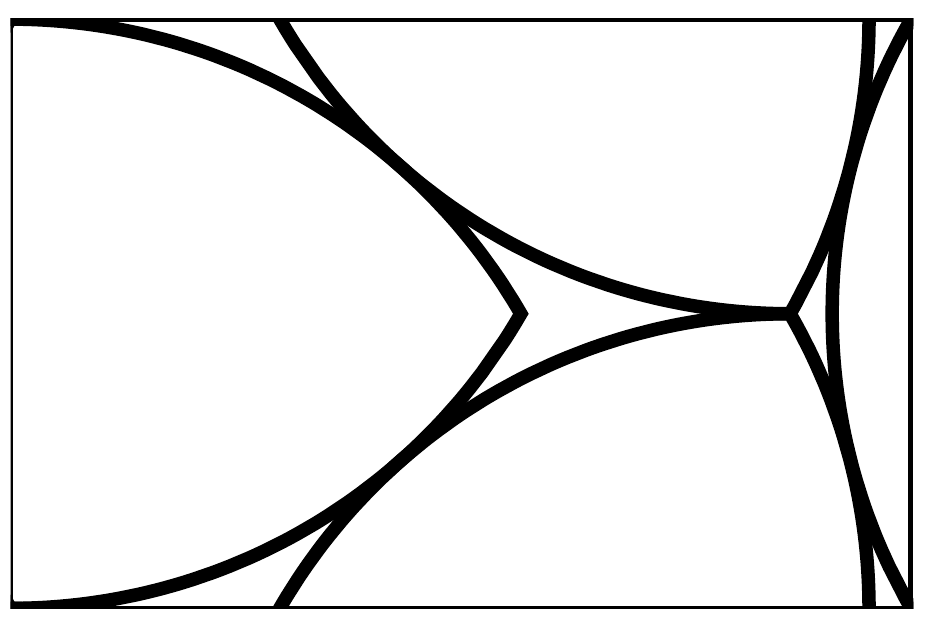}
\caption{Fundamental domain for packing Reuleaux curved triangles. The packing density is $\frac{2 (\pi - \sqrt{3})}{\sqrt{15} + \sqrt{7} - 2 \sqrt{3}} \simeq 0.92288$.}
\label{R3 F domain}
\end{center}
\end{figure}

\begin{figure}[h]
\begin{center}
\includegraphics[width=1.5in]{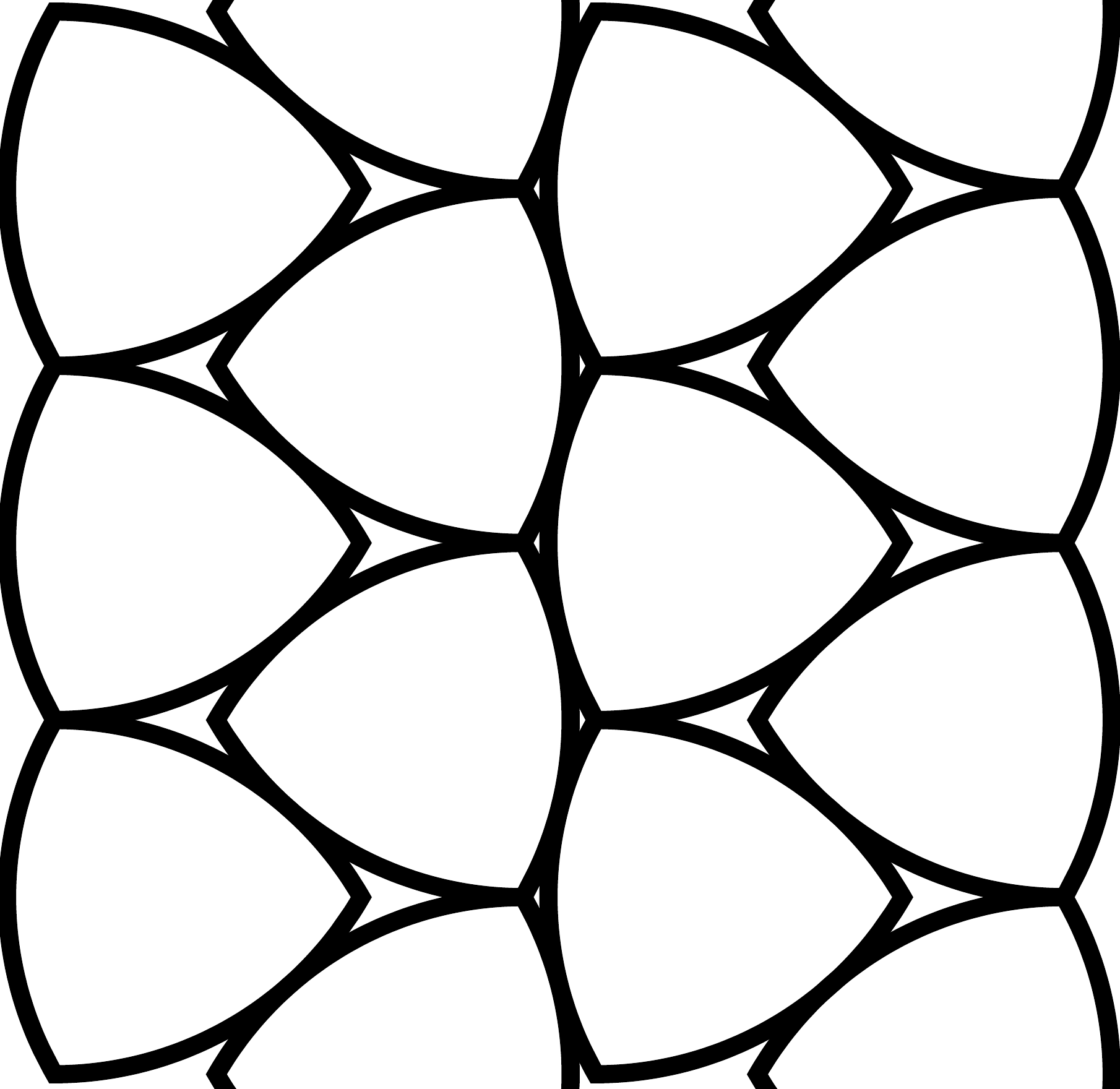}
\caption{A section of a  Reuleaux curved triangle packing. The packing density is $\frac{2 (\pi - \sqrt{3})}{\sqrt{15} + \sqrt{7} - 2 \sqrt{3}} \simeq 0.92288$.}
\label{R3 packing}
\end{center}
\end{figure}


\begin{figure}[h]
\begin{center}
\includegraphics[width=1.8in]{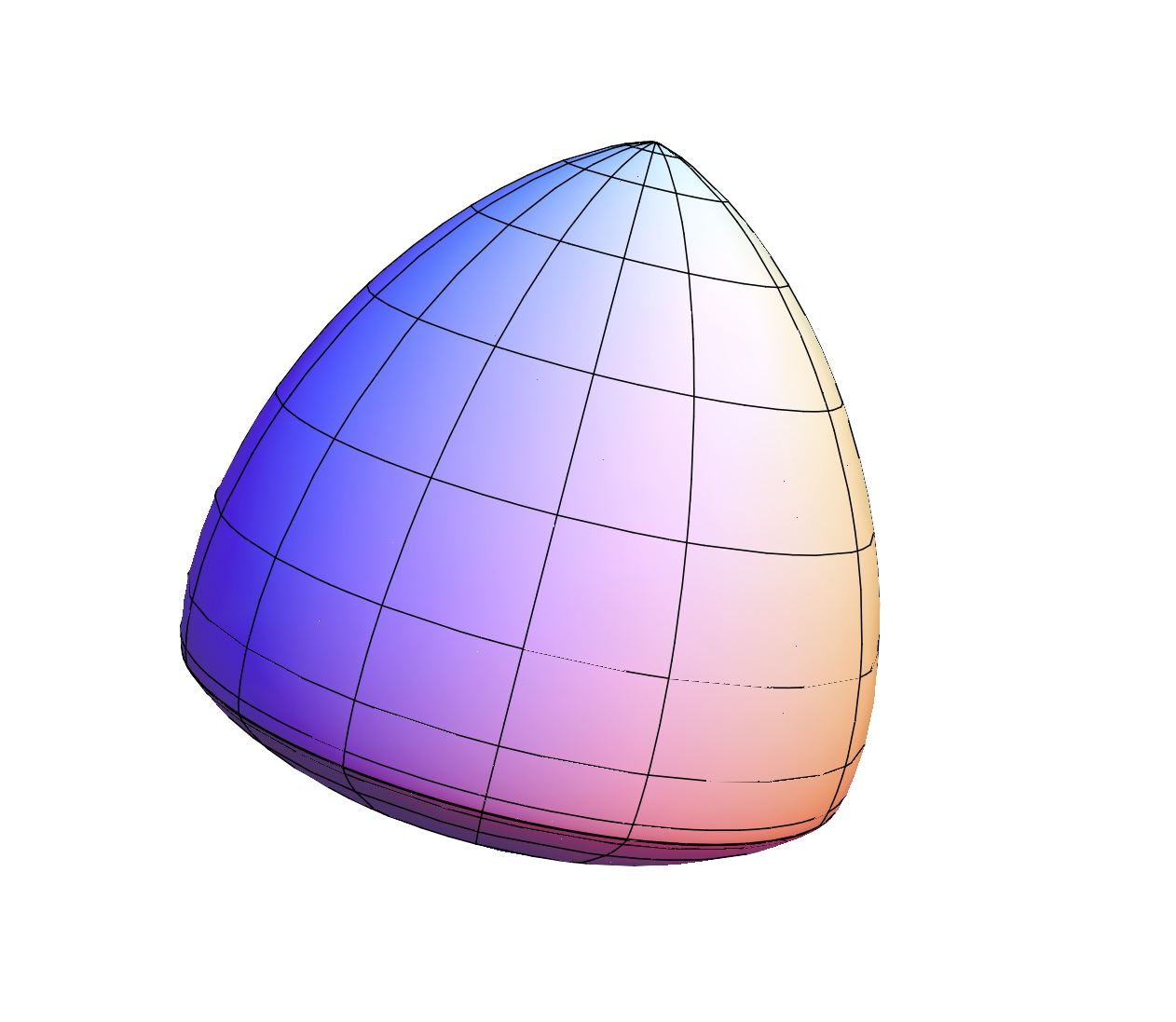}
\caption{Surface of revolution of constant width defined by $p(t,u) = 1 + \frac{1}{8} \cos 3u$.}
\label{HLR_S0}
\end{center}
\end{figure}

\begin{figure}[h]
\begin{center}
\includegraphics[width=4.6in]{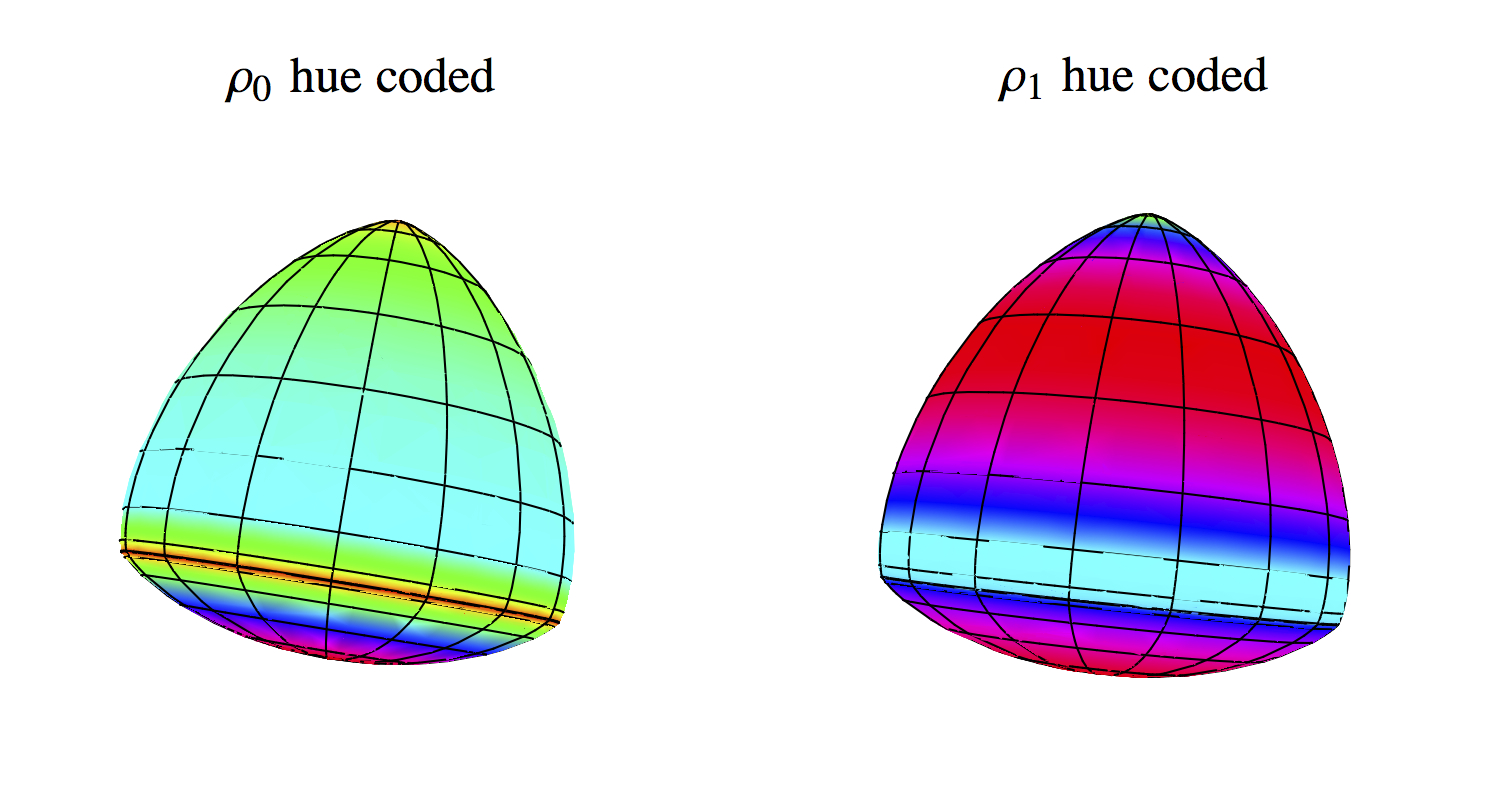}
\caption{Surface of revolution of constant width defined by $p(t,u) = 1 + \frac{1}{8} \cos 3u$. Radii of curvature colored by hue: $0< {\rm red} < {\rm cyan}=1  <  {\rm blue} <2$.}
\label{HLR_S0_hue}
\end{center}
\end{figure}

\begin{figure}[h]
\begin{center}
\includegraphics[width=2.5in]{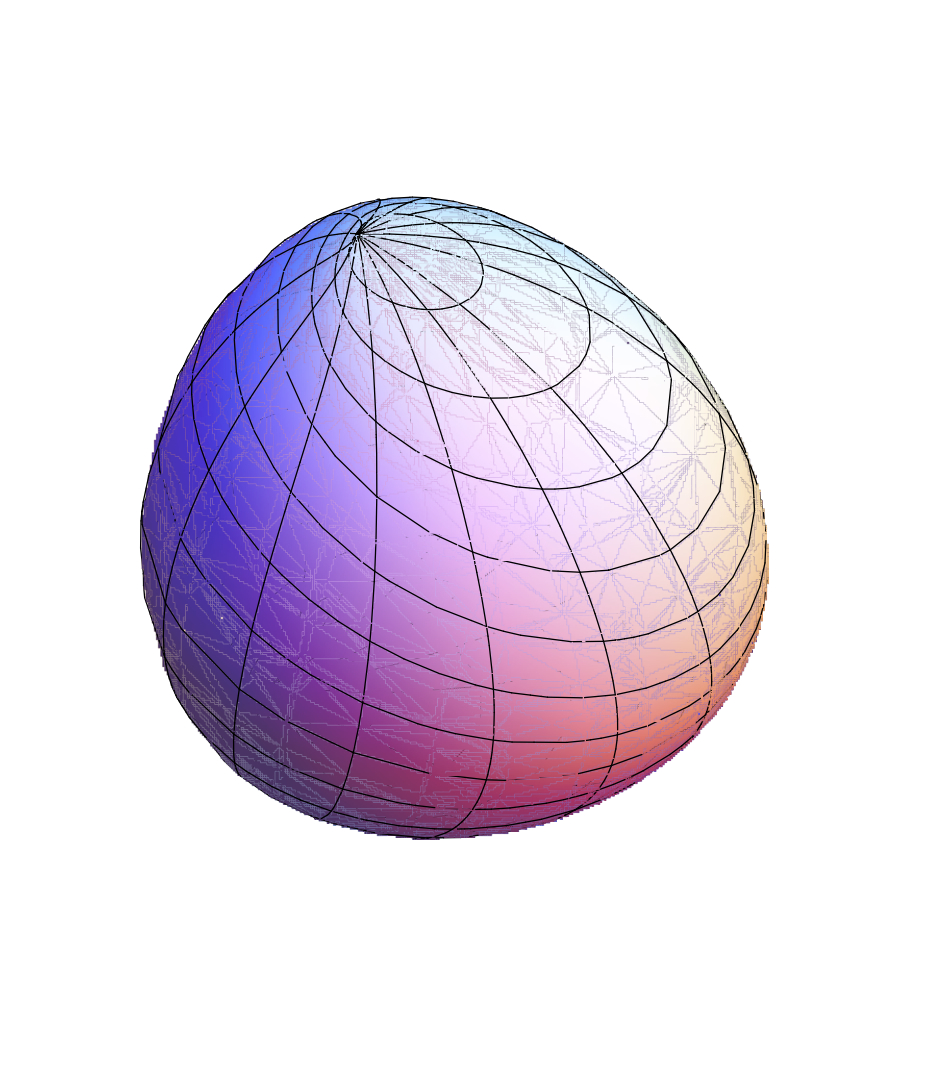}
\caption{Surface of constant width defined by $p(t,u)  = 1 + \frac{1}{2} \cos t \sin^2 u$.}
\label{HLR_S(1,0)}
\end{center}
\end{figure}

\begin{figure}[h]
\begin{center}
\includegraphics[width=2.8in]{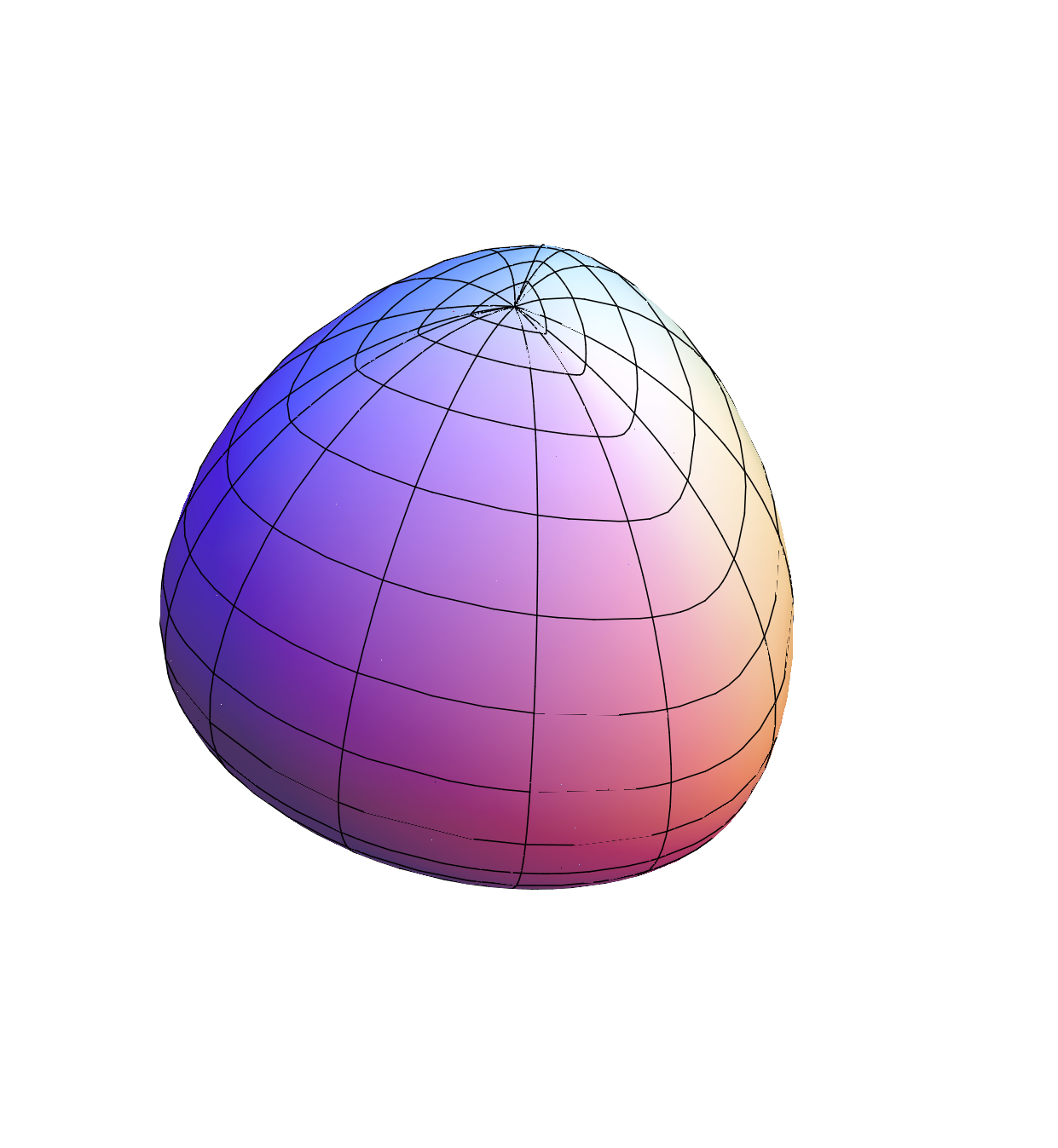}
\caption{Surface of constant width defined by $ p(t,u) = 1 +  \frac{1}{2} \big( \frac{1}{8}\cos 3u +  \frac{1}{8}  \cos 3 t \,  \sin^2 u \big) $.}
\label{HLR_S1}
\end{center}
\end{figure}

\begin{figure}[h]
\begin{center}
\includegraphics[width=2.7in]{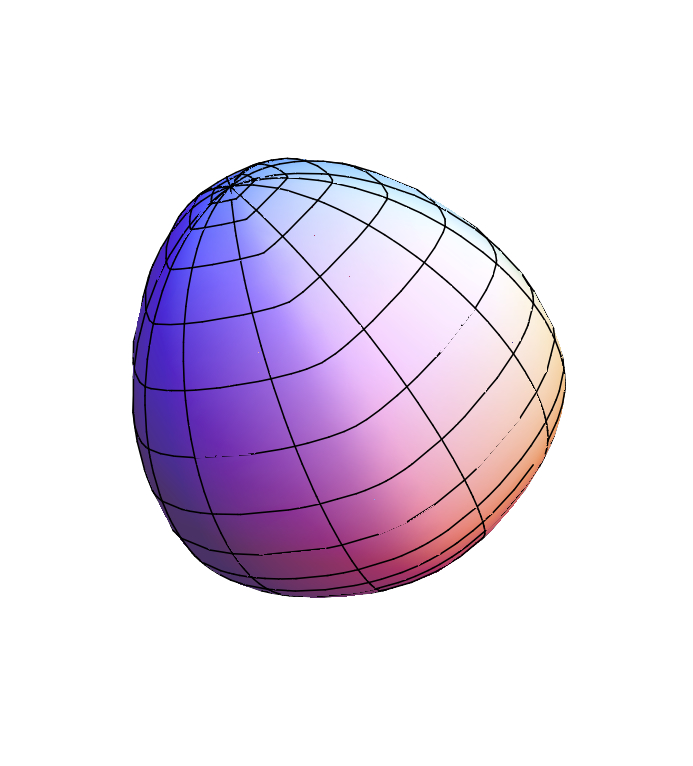}
\caption{Surface of constant width defined by 
 $p(t,u) = 1 + \frac{1}{2} \big( \frac{1}{8} \cos 3 u + \frac{1}{24} \cos 5 t \,  \sin^2 u  \big)$.}
\label{HLR_S(5,3)}
\end{center}
\end{figure}

\begin{figure}[h]
\begin{center}
\includegraphics[width=2.8in]{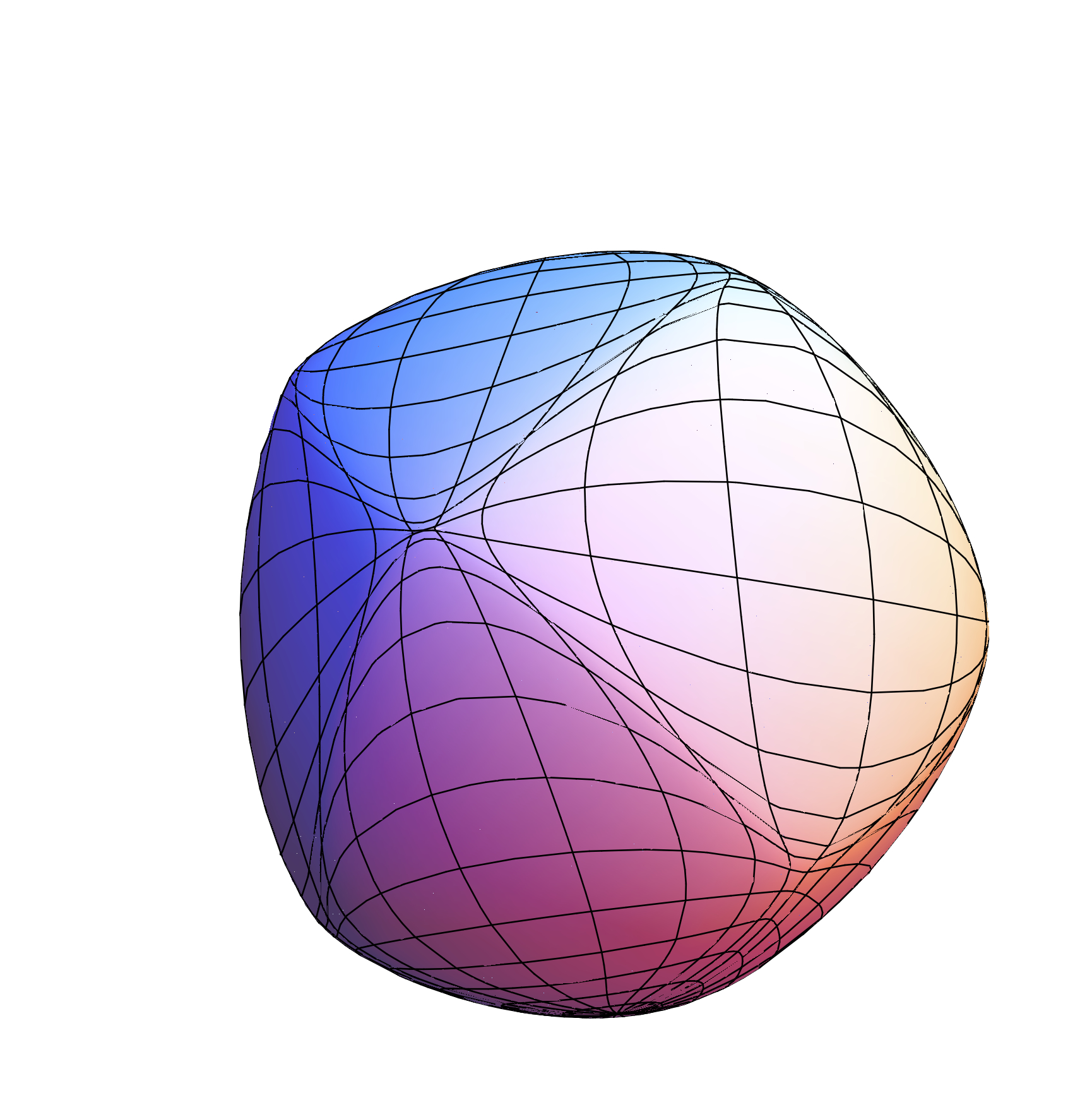}
\caption{Surface of constant width defined by 
 $p(t,u) = 1 +   \frac{1}{10} \cos 3 t \, \sin 3u \,  \sin^2 u $.}
\label{HLR_S(3,3)}
\end{center}
\end{figure}


\end{document}